\documentclass{qjmam}

\usepackage{graphicx}%
\usepackage{amsmath,amssymb,amsthm,eucal,upref,bm}%
\usepackage{labelfig}%

\usepackage{mathrsfs} 

\startpage{1}

\makeatletter
\def\ps@titlepage{%
  \let\@oddhead\@empty
  \let\@evenhead\@oddhead
  \def\@oddfoot{%
      {\vbox{\hsize\textwidth\@parboxrestore\footnotesize
         {\bfseries\journalname}\hfill
         \copyright\ {\bfseries\copyrightname\endgraf}}}%
     \relax
     \relax}%
  \let\@evenfoot\@empty}

\newtheorem{definition}{Definition}[section]

\renewcommand{\geq}{\geqslant} \renewcommand{\leq}{\leqslant}

\DeclareMathAlphabet{\varmathbb}{U}{pxsyb}{m}{n}

\newtheorem{theorem}{Theorem}%

\newcommand{\RR}{\mathbb{R}}%
\newcommand{\D}{\mathrm{d}\kern0.2pt}%
\newcommand{\ii}{\kern0.05em\mathrm{i}\kern0.05em}%
\newcommand{\E}{\textrm{e}}%

\begin{document}

\title[{Babenko's~equation~for~waves~on~water~of~finite~depth}]{Babenko's equation
for periodic gravity waves on water of finite depth: derivation and numerical
solution}

\author[n.~kuznetsov, e.~dinvay] {N. KUZNETSOV \and E. DINVAY}

\address{Laboratory for Mathematical Modelling of Wave Phenomena, \\ Institute for
Problems in Mechanical Engineering, Russian Academy of Sciences, \\ V.O., Bol'shoy
pr. 61, {\rm 199178} St Petersburg, Russian Federation and \\ Department of
Mathematics, University of Bergen, All\'egaten 41, N-5020 Bergen}

\received{8 May 2018, Revised ??}

\maketitle

\eqnobysec

\begin{abstract} 
The nonlinear two-dimensional problem, describing periodic steady waves on water of
finite depth is considered in the absence of surface tension. It is reduced to a
single pseudo-differential operator equation (Babenko's equation), which is
investigated analytically and numerically. This equation has the same form as the
equation for waves on infinitely deep water; the latter had been proposed by Babenko
and studied in detail by Buffoni, Dancer and Toland. Instead of the $2 \pi$-periodic
Hilbert transform $\mathcal{C}$ used in the equation for deep water, the equation
obtained here contains a certain operator $\mathcal{B}_r$, which is the sum of
$\mathcal{C}$ and a compact operator whose dependence on the parameter involves on
the depth of water. Numerical computations are based on an equivalent form of
Babenko's equation derived by virtue of the spectral decomposition of the operator
$\mathcal{B}_r \D / \D t$. Bifurcation curves and wave profiles of the extreme form
are obtained numerically.
\end{abstract}

\section{Introduction} 

Near the end of his remarkable career in both pure and applied mathematics (see
\cite{O} for highlights of major achievements), Konstantin Ivanovich Babenko
(1919--1987) published four brief notes \cite{B,BB,BP,BPR} (the last two appeared
posthumously) on a classical nonlinear problem in the mathematical theory of water
waves, namely, the two-dimensional problem of steady, periodic waves on infinitely
deep water. In this paper dedicated to the centenary of Babenko's birth, we extend
the approach developed in \cite{B} to the case of water of finite depth and deduce a
single pseudo-differential operator equation (Babenko's equation) equivalent to the
corresponding free-boundary problem in some sense explained below (see Section~3.3).
Moreover, using the spectral decomposition of a linear operator involved in the
equation, we transform it to a form convenient for discretisation and then apply a
very robust numerical method that allows us to produce convincing results concerning
global bifurcation branches, secondary bifurcations and free surface profiles
including those of the extreme form.

It was Stokes \cite{S}, who had initiated studies in this field. On
the basis of approximations developed for waves with a single crest per wavelength
(now, they are referred to as {\it Stokes waves}), he made some conjectures about
the behaviour of such waves on deep water. To a great extent, these conjectures had
determined the development of research in the 20th century; see the paper \cite{PT}
and references cited therein to get an idea how these conjectures were proved. In
particular, the so-called Nekrasov's equation was essential for this purpose. The
first version of this nonlinear integral equation for waves on deep water was
derived in \cite{N1} (see also \cite{N3}, Part~1). Soon after that, Nekrasov
generalized his equation to the case of finite depth (see \cite{N2} and \cite{N3},
Part~2). Much later, Amick and Toland \cite{AT2} proposed and investigated a more
sophisticated version of the latter equation.

Levi-Civita \cite{LC} and Struik \cite{St} considered (independently of Nekrasov)
the problem of periodic waves on deep water and on water of finite depth
respectively. The hodograph transform allowed them to reduce the question of
existence of waves to that of finding a pair of conjugate harmonic functions
satisfying nonlinear Neumann boundary conditions. The existence proofs given in
\cite{LC} and \cite{St} are based on a majorant method for demonstrating the
convergence of power series that provide formal solutions. In his book \cite{Z2},
Chapter~71, Zeidler writes about these proofs that they are `very complicated' in
view of `voluminous computations' involved. By now, both techniques (Nekrasov's
equations and the method of Levi-Civita and Struik) are investigated in detail by
means of analytic bifurcation theory. An account of this theory can be found in the
books \cite{Z1} and \cite{BT}, whereas many authors studied its application to
equations describing steady water waves; these results are summarised in \cite{T}
(deep water) and in \cite{Z2}, Chapter~71 (water of finite depth), where one also
finds detailed historical remarks. It should be mentioned that Krasovskii
\cite{Kras} extended studies of water waves to the case of a periodic wavy bottom.

Another method for periodic waves on deep water (with and without surface tension)
was developed in detail by Buffoni, Dancer and Toland \cite{BDT1,BDT2}. In the
absence of surface tension, it is based on the so-called Babenko's equation:
\begin{equation}
\mu \, \mathcal{C} (v') = v + v \, \mathcal{C} (v') + \mathcal{C} (v' v) , \quad t \in
(-\pi, \pi) . \label{bid}
\end{equation}
Here $\mu$ is a dimensionless bifurcation parameter (it is related to the speed of
wave pro\-pagation), which must be found along with the $2 \pi$-periodic function
$v (t)$ that describes the wave profile parametrically in certain dimensionless
coordinates; $'$ stands for differentiation with respect to $t$ and ${\cal C}$ is
the $2 \pi$-periodic Hilbert transform (the conjugation operator in the theory of
Fourier series; see, for example, \cite{Z}). It is defined on $L^2 (-\pi, \pi)$ by
linearity from the following relations:
\begin{equation}
{\cal C} (\cos n t) = \sin n t \ \ \mbox{for} \ n \geq 0 , \quad {\cal C} (\sin n t)
= - \cos n t \ \ \mbox{for} \ n \geq 1 . \label{HT}
\end{equation}

The original form of \eqref{bid} was announced by Babenko \cite{B} (see also
\cite{OS}, Section~3.7), where the equation is derived and expressed in terms of the
self-adjoint operator $J_0 = \mathcal{C} \, \D / \D t$. In his second note
\cite{BB}, Babenko outlines how to prove the local existence theorem for his
equation in a neighbourhood of the first bifurcation point equal to unity. He
demonstrates that $\mu$ must be equal to $1 + \epsilon^2$, and, after changing the
unknown function by applying the invertible operator $I + J_0$ to the original one,
a solution is sought in the form of expansion in powers of $\epsilon$. It is shown
that the expansion converges for $\epsilon \leq 1 / 25$. Some numerical computations
related to the Babenko's version of equation \eqref{bid} are presented in
\cite{BP,BPR}.

A generalization of Babenko's equation was later studied in \cite{ST}. Besides,
Longuet-Higgins \cite{LH1} derived an infinite system of algebraic equations
equivalent to \eqref{bid} (see also \cite{BS} and \cite{OS}, Section~3.6). He used
this system for numerical computations of Stokes-wave bifurcations (see \cite{LH2}
and also \cite{OS}, Sections~3.2 and 3.8). It is worth mentioning that this system
naturally appears from the Lagrangian formalism developed in \cite{Ba}. In the paper
\cite{LH3}, a similar system of quadratic relations between the Fourier coefficients
of the wave profile was obtained in the case of water having finite depth.

Interesting results concerning the secondary or sub-harmonic bifurcations from
branches describing Stokes-wave solutions of \eqref{bid} are proved in the articles
\cite{BDT1} and \cite{BDT2}. In the first of these, it is shown that such
bifurcations do not occur in a neighbourhood of those points, where Stokes waves
bifurcate from a trivial solution. On the other hand, significant numerical evidence
about the existence of steady periodic waves that distinguish from Stokes waves had
appeared in the 1980s. These new waves have more than one crest per period and
bifurcate from Stokes waves. Branches of sub-harmonic bifurcations for deep water
were first computed by Chen and Saffman \cite{CS}, whereas Vanden-Broeck \cite{VB}
obtained similar result for water of finite depth by solving numerically an
integro-differential system arising after the hodograph transform; this system was
proposed in \cite{VBS}. Craig and Nicholls \cite{CN} used a different method for
computing numerical results generalising those of Vanden-Broeck. Moreover, it was
shown that non-symmetric periodic waves exist on water of finite depth for which
purpose a weakly nonlinear Hamiltonian model was applied in \cite{Zuf1}.

References to other works containing numerical results on sub-harmonic bifurcations
can be found in \cite{OS} and \cite{BM}. In the latter paper, some theoretical
insights concerning these bifurcations are also given. The above mentioned numerical
observations were confirmed rigorously in \cite{BDT2}, where it was `concluded that
the sub-harmonic bifurcations [\dots] are {\it an inevitable consequence} of the
formation of Stokes highest wave'. A characteristic property of the latter wave is
the angle equal to $2\pi/3$ formed at the crest by two smooth, symmetric curves.
Concerning this wave see \cite{PT} and references cited therein.

Apart of numerical approaches mentioned above, the direction of studies was quite
different for water of finite depth. Of course, Nekrasov's equation and the approach
of Struik were both developed for waves of a given wavelength. On the other hand, a
psuedo-differential equation in terms of variables arising after the hodograph
transform was derived in \cite{KK}. This equation describes all waves for which the
rate of flow per unit span and the Bernoulli constant are prescribed and serves for
justifying the Benjamin--Lighthill conjecture for near-critical waves. However, it
is not suitable for studying the Stokes-wave and sub-harmonic bifurcations. The
results obtained for equation \eqref{bid} in \cite{BDT1,BDT2} show that Babenko's
equation serves better for this purpose. Here the analysis of equation similar to
\eqref{bid}, but describing waves on water of finite depth is initiated and new
numerical results are obtained.

Besides, Constantin, Strauss and V\u{a}rv\u{a}ruc\u{a} investigated the problem of
water waves with constant vorticity in their recent paper \cite{CSV}. In the case of
finite depth, this problem is reduced to a quasilinear pseudo-differential system
provided the vorticity is non-zero. In the irrotational case (that is, for zero
vorticity) and for some particular value of a parameter involved in the system, the
latter turns into a single equation that has the same form as \eqref{bid} with
$\mathcal C$ changed to the so-called periodic Hilbert transform on a strip. In
Section~5, we compare this equation with Babenko's equation derived in this paper;
see \eqref{37} below.

The plan of the paper is as follows. For the problem describing steady, periodic
waves on water of finite depth two equivalent statements\,---\,dimensional and
non-dimensional\,---\,are formulated in Section~2. Babenko's equation is derived
from the non-dimensional for\-mulation in Section~3.1 and the existence of local
bifurcation branches for this equation is proved on the basis of the
Crandall--Rabinowitz theorem in Section~3.2. In Section~3.3, we outline how to
obtain a solution of the non-dimensional problem from a given solution of Babenko's
equation. Numerical procedure applied for solving Babenko's equation is presented in
Section~4 along with various bifurcation curves and wave profiles obtained with its
help. Section~5 contains concluding remarks.

\section{Statements of the problem}

In its simplest form, the problem of steady surface waves concerns the
two-dimensional, irrotational motion of an inviscid, incompressible, heavy fluid,
say water, bounded above by a free surface and below by a rigid horizontal bottom.
(For example, this kind of motion occurs in water occupying an infinitely long
channel with rectangular cross-section and having uniform width.) In an appropriate
frame of reference the velocity field of steady motion is time-independent as well
as the free-surface profile, and they can be described in two equivalent ways that
differ by prescribed parameters.

\subsection{The Benjamin--Lighthill statement for steady waves}

The classical formulation proposed by Benjamin and Lighthill is convenient for
justification of their conjecture (see \cite{Ben} and \cite{KK}, where it had been
justified for Stokes waves and all near-critical waves respectively). In this
formulation, $Q$\,---\,the constant rate of flow per channel's unit span\,---\,is
prescribed along with the total head $R$ also referred to as the Bernoulli constant.
Let Cartesian coordinates $(X,Y)$ be chosen so that the bottom coincides with the
$X$-axis and gravity acts in the negative $Y$-direction, whereas the wave profile
has a crest on the $Y$-axis (this does not restrict generality). Thus, the profile
is given by the graph of an unknown positive function $\xi$ (that is, $Y = \xi (X)$,
$X \in \RR$), attaining its maximum at $X=0$. Moreover, we suppose that $\xi$ is
continuously differentiable and even. In the longitudinal section of the water
domain ${\cal D} = \{ X \in \RR,\ 0 < Y < \xi (X) \}$, the velocity field is
described by the stream function $\Psi (X, Y)$, that is, the projections of the
velocity vector at $(X, Y)$ on the $X$- and $Y$-axis are $\Psi_Y$ and $-\Psi_X$
respectively. It is assumed that $\Psi$ belongs to $C^2 ({\cal D}) \cap C^1
(\bar{\cal D})$ and is an even function of $X$ (hence it is bounded on $\bar{\cal
D}$).

If the surface tension is neglected, then $\Psi$ and $\xi$ must satisfy the
following free-boundary problem:
\begin{eqnarray}
&& \Psi_{XX} + \Psi_{YY} = 0, \quad (X,Y) \in {\cal D}; \label{1} \\ && \Psi (X, 0)
= - Q, \quad X \in \RR; \label{2} \\ && \Psi (X, \xi (X)) = 0, \quad X \in \RR;
\label{3} \\ && \frac{1}{2} |\nabla \Psi (X, \xi (X))|^2 + g \, \xi (X) = R , \quad  
X \in \RR . \label{4}
\end{eqnarray}
In the left-hand side of the last relation usually referred to as Bernoulli's
equation, $g > 0$ is the constant acceleration due to gravity. It is known that
non-trivial solutions of problem (\ref{1})--(\ref{4}) exist only when $Q \neq 0$ and
$R > R_c = \frac{3}{2} (Q g)^{2/3}$. For the proof of the first relation see
Proposition 1.1 in \cite{KK1}, whereas the last inequality is proved in \cite{KK2}
under weaker assumptions than listed above. In what follows, these restrictions on
$Q$ and $R$ hold; moreover, we suppose (without loss of generality) that $Q > 0$.

\subsection{A non-dimensional statement for periodic waves}

Let us assume that $\xi$ is a $2 \ell$-periodic function ($\ell > 0$), whereas $\Psi
(X, Y)$ is $2 \ell$-periodic in $X$, but the constant $R$ is to be found along with
these functions from problem (\ref{1})--(\ref{4}). In order to reduce the
reformulated problem to a non-dimensional form, we average Bernoulli's equation over
$(-\ell, \ell)$. Since $\Psi$ is constant on the free surface, we obtain that $c^2 =
2 (R - g H)$, where
\begin{equation}
H = \frac{1}{2 \ell} \int_{-\ell}^\ell \xi (X) \, \D X \quad \mbox{and} \quad c^2 =
\frac{1}{2 \ell} \int_{-\ell}^\ell \left| \frac{ \partial \Psi}{\partial n} (X, \xi
(X)) \right|^2 \D X . \label{abe}
\end{equation}
Here $n$ is the unit normal to $Y = \xi (X)$ directed out of ${\cal D}$. One can
show that the last equality \eqref{abe} is true with the same constant $c^2$ when
this curve is changed to $Y = \tilde \xi (X)$\,---\,an arbitrary
streamline\,---\,and $n$ is changed to $\tilde n$\,---\,the unit normal to this
streamline. Therefore, $c > 0$ is the unknown mean velocity of flow in the positive
direction of the $X$-axis.

It is convenient to introduce the following non-dimensional quantities: $h = \pi H /
\ell$ (the mean depth of flow) and $Q_0 = Q / \sqrt{ g (\ell / \pi)^3}$ (the rate
of flow). Now we scale the dimensional variables and shift the vertical variables
downwards as follows:
\begin{equation}
x = \frac{\pi}{\ell} X ,\ y = \frac{\pi}{\ell} Y - h ; \quad \eta (x) = \frac{\pi}
{\ell} \, \xi (X) - h ; \quad \psi (x,y) = \frac{Q_0}{Q} \Psi (X,Y) . \label{dlv}
\end{equation}
This is advantageous because the new unknown $\eta$ is a $2 \pi$-periodic and even
function of $x$ satisfying the following condition:
\begin{equation}
\int_{-\pi}^\pi \eta (x) \, \D x = 0 . \label{eta}
\end{equation}
Furthermore, the function $\psi$ has the same properties on $\bar D$ as $\Psi$ has
on $\bar{\cal D}$; namely,
\[ \psi \in C^1 (\bar{D}) \cap C^2 (D), \quad \mbox{where} \
D = \{ x \in \RR, -h < y < \eta (x) \},
\]
and is a $2 \pi$-periodic and even function of $x$. Moreover, the change of
variables (\ref{dlv}) reduces relations (\ref{1})--(\ref{4}) to the following
\begin{eqnarray}
&& \psi_{xx} + \psi_{yy} = 0, \quad (x,y) \in D; \label{lapp} \\ && \psi (x, -h) =
-Q_0 , \quad x \in \RR; \label{bcp} \\ && \psi (x, \eta (x)) = 0, \quad x \in \RR;
\label{kcp} \\ && |\nabla \psi (x, \eta (x))|^2 + 2 \eta (x) = \mu , \quad x \in \RR . 
\label{bep}
\end{eqnarray}
In the non-dimensional Bernoulli equation, the parameter $\mu = \pi c^2 / (g \ell)$
is the Froude number squared which must be found along with $\eta$ and $\psi$.
Besides, $\mu / 2$ serves as the independent of $h$ upper bound for $\eta$ with
equality achieved only for the wave of extreme form with the Lipschitz crest; see
\cite{CN}. Thus, the non-dimensional statement of the problem is as follows.

\begin{definition}
{\rm Let $Q_0$ and $h$ be given positive numbers, then problem P$(Q_0,h)$ is to find
a triple $(\mu, \eta, \psi)$ from relations (\ref{lapp})--(\ref{bep}) so that $\mu$
is positive, $\eta$ satisfies condition \eqref{eta}, whereas all other properties of
$\eta$ and $\psi$ (smoothness, $2 \pi$-periodicity and symmetry) are as described above.}
\end{definition}

On the other hand, having a solution of problem P$(Q_0,h)$, formulae (\ref{dlv})
yield a $2 \ell$-periodic solution of problem (\ref{1})--(\ref{4}) for any $\ell >
0$, for which purpose one has to put $Q = Q_0 \sqrt{ g (\ell / \pi)^3}$ and to
determine $R$ from the equality $c^2 = 2 (R - g H)$ with $c^2 = \mu g \ell / \pi$
and $H = h \ell / \pi$.

\section{Babenko's equation}

The aim of this section is to derive a single nonlinear pseudo-differential equation
that has the same form as (\ref{bid}), but $\mathcal{C}$ is replaced by the sum of
$\mathcal{C}$ and some compact operator depending on a real parameter. The equation
is related to the family of problems P$(Q_0,h)$ in the following sense. The value of
operator's parameter together with equation's solution allow us to determine $h$ and
to obtain some solution of the water-wave problem.

\subsection{Derivation of Babenko's equation}

First, we follow considerations in Section~8 of \cite{N3}; see also the rather
recent paper \cite{Bod}. By $w (z) = \varphi + \ii \psi$ $(z = x + \ii y)$ we denote
the complex potential restricted to the one-wave domain $D_{2 \pi} = \{ -\pi < x <
\pi , -h < y < \eta (x) \}$ of some periodic wave on water of a certain depth $h$.
Here $\varphi (x, y)$ is the odd in $x$ harmonic conjugate to $\psi$, for which
purpose the additive constant must be chosen properly. For some $r \in (0,1)$ we
consider a conformal mapping, say $u (z)$, from the $z$-plane to the auxiliary
$u$-plane; it maps $D_{2 \pi}$ onto
\begin{equation}
A_r = \{ r < |u| < 1 ; \ \Re\,u \notin (-1 , -r) \ \mbox{when} \ \Im\,u =0 \} .
\label{A_r}
\end{equation} 
This annular domain has a cut which makes it simply connected and the map is such
that the images of the upper and bottom parts of $\partial D_{2 \pi}$ are
\[ \{ |u|=1 ; \Re\,u \neq -1 \} \quad \mbox{and} \quad \{ |u|=r ; \Re\,u \neq -r \}
\]
respectively, whereas the left (right) side of $\partial D_{2 \pi}$ is mapped onto
the upper (lower respectively) side of the cut $\{ \Re\,u \in (-1 , -r) \
\mbox{when} \ \Im\,u =0 \}$. Putting
\begin{equation}
u = \E^{-\ii w} \quad \mbox{and} \quad \frac{\D z}{\D u} = \ii \left[ u^{-1} + f (u)
\right] \, , \label{zu}
\end{equation}
where $f (u)$ is a certain Laurent series, one obtains that
\begin{equation}
- \frac{\D w}{\D z} = \left[ 1 + u f (u) \right]^{-1} \, . \label{wu}
\end{equation}
In \cite{N3}, Section 8, this formula serves as the basis for deriving Nekrasov's
equation in the case of finite depth. An equivalent form of this equation is given
in \cite{Bod}; see equation (1.1) there.

According to the second equality \eqref{zu}, the general form of the inverse
conformal map\-ping $A_r \ni u \mapsto z \in D_{2 \pi}$ is as follows:
\begin{equation}
z (u) = \ii \Big[ \log u - a_0 + \sum_{k=1}^\infty a_k \left( u^k - r^{2 k} u^{-k}
\right) \Big] \, , \quad \mbox{where} \ a_k \in \RR , \ k = 0,1,2,\dots \, .
\label{z_u}
\end{equation}
Here, we put minus at $a_0$ because it is convenient in what follows. The fact that
all coefficients $a_k$ are real is a consequence of equality \eqref{wu} because
$\psi$ is equal to a real constant on the bottom part of $\partial D_{2 \pi}$ which
corresponds to $\{ |u|=r ; \Re\,u \neq -r \}$\,---\,the circumference cut on the
negative real axis.

Let us derive some relations for the coefficients from \eqref{z_u}. First, for $u=r$
we obtain
\begin{equation}
a_0 = h + \log r . \label{a_0}
\end{equation}
Substituting $u = \E^{\ii t}$, $t \in (- \pi, \pi)$, into \eqref{z_u} and separating
the real and imaginary parts, we arrive at the following parametric representation
of the free surface profile:
\begin{equation}
x (t) = - t - \sum_{k=1}^\infty a_k \left( 1 + r^{2 k} \right) \sin kt \, , \quad
\eta (t) = - a_0 + \sum_{k=1}^\infty a_k \left( 1 - r^{2 k} \right) \cos kt \, .
\label{x_eta_t}
\end{equation}
Now we see that another relation is equivalent to condition \eqref{eta} written in
terms of the last two expressions, namely:
\begin{equation}
\int_{-\pi}^\pi \!\! \eta (t) x' (t) \, \D t = 0 \quad \Longleftrightarrow \quad a_0 =
\frac{1}{2} \sum_{k=1}^\infty k \, a_k^2 \left( 1 - r^{4 k} \right) \, .
\label{etat}
\end{equation}
It follows from the last equality that $a_0 > 0$ in the non-trivial case. Then
equality \eqref{a_0} shows that the value of $r$ is related not only to the depth
$h$, but also to a particular solution of problem P$(Q_0,h)$.

It is worth mentioning that both expressions \eqref{x_eta_t} are similar to those
for the infinite depth (cf.~\cite{OS}, Section~3.7, where Babenko's results are
outlined), and in that case, a consequence is the relation $x_t = - 1 - \mathcal{C}
\eta_t$ with
\[ (\mathcal{C} v) (t) = \frac{1}{2 \pi} \int_{-\pi}^\pi v (\tau) \cot \frac{t-\tau}{2}
\D \tau \, ,
\] 
which is the form of the $2 \pi$-periodic Hilbert transform alternative to formulae
\eqref{HT}.

The crucial point for obtaining a similar relation in the case of finite depth is to
introduce the operator $\mathcal{B}_r = \mathcal{C} + \mathcal{K}_r$ for $r \in (0,
1)$, where
\begin{equation}
(\mathcal{K}_r v) (t) = \frac{2}{\pi} \int_{-\pi}^\pi v (\tau) K_r (t-\tau) \, \D
\tau  \quad \mbox{with} \ \ K_r (t-\tau) = \sum_{n=1}^\infty \frac{r^{2 n}}{1 - r^{2
n}} \sin  (t-\tau) . \label{K_r}
\end{equation}
It is straightforward to calculate that $\mathcal{B}_r$ can also be defined on $L^2
(-\pi, \pi)$ by linearity from the following relations
\begin{equation}
\mathcal{B}_r (\cos n t) = \frac{1 + r^{2 n}}{1 - r^{2 k}} \sin n t \ \ \mbox{for} \
n \geq 0 , \quad \mathcal{B}_r (\sin n t) = - \frac{1 + r^{2 n}}{1 - r^{2 n}} \cos n
t \ \ \mbox{for} \ n \geq 1 \label{HTB_r}
\end{equation}
that are similar to \eqref{HT}. Combining these formulae and \eqref{x_eta_t} yields
that
\begin{equation}
x_t = - 1 - \mathcal{B}_r \eta_t \quad \mbox{for} \ t \in (- \pi, \pi) \, .
\label{30}
\end{equation}

An important fact about the operator $\mathcal{B}_r$ is that it is a conjugation in
the following sense. If $F (u)$ is analytic in $A_r$ and $\Im F$ vanishes
identically on $\{ |u|=r ; \Re\,u \neq -r \}$, then
\begin{equation}
\Re F (\E^{\ii t}) + [ \mathcal{B}_r (\Im F) ] (t) = 0 \quad \mbox{for all} \ t
\in (- \pi, \pi) . \label{32}
\end{equation}

Let us calculate the derivative $z_\varphi$ of the mapping inverse to the complex
potential. In view of the first equality \eqref{zu}, we have
\[ z_\varphi = z_u \, u_\varphi = - \ii z_u \, \E^{-\ii w} w_\varphi = - \ii u \, z_u .
\]
Combining this and \eqref{z_u}, we obtain that
\begin{equation}
z_\varphi = 1 + \sum_{k=1}^\infty k a_k \left( u^k + r^{2 k} u^{-k} \right) \, ,
\label{33}
\end{equation}
and the function on the right-hand side is analytic in $A_r$. Since $z_\varphi$ does
not vanish in the closure of $A_r$, we have that $z_\varphi^{-1} = |\nabla
\varphi|^2 \, \overline{z_\varphi}$ is also analytic in $A_r$. Moreover, the
Bernoulli equation \eqref{bep} implies that
\begin{equation}
z_\varphi^{-1} = (\mu - 2 \eta) (x_\varphi - \ii y_\varphi) = (\mu - 2 \eta) (\ii
\eta_t - x_t) \quad \mbox{when} \ u = \E^{-\ii t} ,
\label{34}
\end{equation}
(cf.~formula (3.38) in \cite{OS}). Here the second equality is a consequence of the
Cauchy--Riemann equations. Then equality \eqref{30} yields that
\begin{equation}
z_\varphi^{-1} = (\mu - 2 \eta) (1 + \mathcal{B}_r \eta_t + \ii \eta_t) \quad
\mbox{for} \ t \in (-\pi, \pi) .
\label{36}
\end{equation}

It follows from previous considerations that the constant in the Laurent expansion
of $z_\varphi^{-1}$ is equal to $\mu$. Furthermore, $\Im \{z_\varphi^{-1} - \mu\}$
vanishes identically on $\{ |u|=r ; \Re\,u \neq -r \}$, which allows us to apply
formula \eqref{32} to the function $z_\varphi^{-1} - \mu$, whose trace on $\{u =
\E^{\ii t}\}$ is equal to
\[ (\mu - 2 \eta) \, \mathcal{B}_r (\eta') - 2 v + \ii (\mu - 2 \eta) \eta' .
\]
Here again $'$ stands for differentiation with respect to $t$. Thus, we arrive at
\[ (\mu - 2 \eta) \, \mathcal{B}_r (\eta') - 2 \eta + \mathcal{B}_r \, [(\mu - 2 \eta)
\eta'] = 0 \quad \mbox{for} \ t \in (-\pi, \pi)\, ,
\]
which simplifies to Babenko's equation for waves on water of finite depth:
\begin{equation}
\mu \, \mathcal{B}_r (\eta') = \eta + \eta \, \mathcal{B}_r (\eta') + \mathcal{B}_r \,
(\eta' \eta) \quad \mbox{for} \ t \in (-\pi, \pi) .
\label{37}
\end{equation}

This equation is similar to \eqref{bid} and the derivation procedure yields that for
each $r \in (0, 1)$ it is related to some solution of problem P$(Q_0,h)$.

\subsection{Local bifurcation branches for Babenko's equation}

To show the existence of small solutions of equation \eqref{37}, bifurcating from
the zero solution, we apply the Crandall--Rabinowitz theorem (see \cite{CR},
Theorem~1.7) that deals with bifurcation from simple eigenvalues of the linearised
equation; its formulation is as follows.

\begin{theorem} 
Let ${\cal X}$, ${\cal Y}$ be real Banach spaces with the continuous embedding
${\cal X} \subset {\cal Y}$. If a continuous map ${\cal F} (\mu, v): \RR \times
{\cal X} \mapsto {\cal Y}$ has the following properties:

\vspace{2mm}

{\rm (i)} the equality ${\cal F} (\mu, 0) = 0$ holds for all $\mu \in \RR$,

{\rm (ii)} the operators ${\cal F}_\mu$, ${\cal F}_v$ and ${\cal F}_{\mu v}$ exist
and are continuous,

{\rm (iii)} for some $\mu^*$ the operator ${\cal F}_v (\mu^*, 0)$ is a Fredholm one
with zero index and its null-space is one-dimensional,

{\rm (iv)} if the null-space of ${\cal F}_v (\mu^*, 0)$ is generated by $v^{(0)}$,
then ${\cal F}_{\mu v} (\mu^*, 0) \, v^{(0)}$ does not belong to the range of ${\cal
F}_v (\mu^*, 0)$.

\vspace{2mm}

\noindent Then a sufficiently small $\varepsilon > 0$ exists and a continuous curve
\[ \{ (\mu (s), \, v (s)) : |s| < \varepsilon \} \subset \RR \times {\cal X} , 
\]
bifurcates from $(\mu^*, 0);$ for pairs belonging to this curve
\[ \mu (s) = \mu^* + o (s) \quad and \quad v (s) = s \, v^{(0)} + o (s)
\quad when \ 0 < |s| < \varepsilon .
\]
Moreover, if ${\cal F}_{v v}$ is continuous, then the curve is of class $C^1$.
\end{theorem}

As in \cite{BDT1}, we say that a real-valued function $v$ belongs to the Sobolev
space $H_0$ provided it is absolutely continuous on $[-\pi, \pi]$, $v (-\pi) = v
(\pi)$, and its weak derivative $v'$ belongs to $L^2 (-\pi, \pi)$. Let $\hat H_0$ be
the subspace of $H_0$ consisting of even functions.

In terms of the map ${\cal F}: \RR \times \hat H_0 \mapsto L^2 (-\pi, \pi)$ defined
by
\begin{equation}
{\cal F} (\mu, v) = \mu \mathcal{B}_r (v') - v + v \, \mathcal{B}_r (v') +
\mathcal{B}_r \, (v' v) ,
\label{calF}
\end{equation}
equation \eqref{37} takes the following form:
\begin{equation}
{\cal F} (\mu, v) = 0 , \quad (\mu, v) \in \RR \times \hat H_0 .
\label{eqF}
\end{equation}
Let us apply the Crandall--Rabinowitz theorem to this equation to obtain local
branches of Stokes-wave solutions of small amplitude, for which purpose we have to
check conditions (i)--(iv) for ${\cal F} (\mu, v)$.

It is clear that (i) and (ii) are fulfilled and ${\cal F}_{v} (\mu, 0) = \mu
\mathcal{B}_r \, (\D / \D t) - I$, where $I$ is the identity operator. Hence the set
of bifurcation points of equation \eqref{eqF} is $\{ \mu_n \}_{n=1}^\infty$, where
\begin{equation}
\mu_n = \frac{1 - r^{2n}}{n (1 + r^{2n})} , \quad n=1,2,\dots \, ,
\label{lambda_n}
\end{equation}
are the characteristic values of $\mathcal{B}_r \, (\D / \D t)$. Furthermore, ${\cal
F}_{v} (\mu_n, 0)$ is a Fredholm operator, its index is equal to zero for every
$\mu_n$, and the corresponding null-space in $\hat H_0$ is one-dimensional being
generated by $v^{(0)}_n (t) = \cos n t$, thus yielding condition (iii). Since ${\cal
F}_{\mu v} (\mu_n, 0) = - I$ we see that ${\cal F}_{\mu v} (\mu_n, 0) v^{(0)}_n (t)
= - \cos n t$, and so condition (iv) for $n=1,2,\dots$ is a consequence of the fact
that the equation
\[ \mu_n \mathcal{B}_r (v') - v = - \cos n t
\]
has no solution. Indeed, a solution of this equation exists if and only if its
right-hand side is orthogonal to the null-space of the adjoint operator
\[ [\mu_n \mathcal{B}_r \, (\D / \D t) - I]^* = \mu_n (\D / \D t) \mathcal{B}_r - I .
\]
Since its null-space is one-dimensional and generated by $\cos n t$, the
orthogonality condition does not hold for $- \cos n t$. This completes verification
of condition (iv).

Then the Crandall--Rabinowitz theorem yields the following.

\begin{theorem}
For every $n=1,2,\dots$ there exists $\varepsilon_n > 0$ such that for $0 < |s| <
\varepsilon_n$ there is the family $\big( \mu_n^{(s)}, \, v_n^{(s)} \big)$ of
Stokes-wave solutions to equation \eqref{eqF}. Together with the bifurcation point
$(\mu_n , 0)$, where $\mu_n$ is given by formula \eqref{lambda_n}, the points of
this family form the continuous curve
\[ C_n = \big\{ \big( \mu_n^{(s)}, \, v_n^{(s)} (t) \big) : |s| < \varepsilon_n
\big\} \subset \RR \times \hat H_0 , \quad n=1,2,\dots \, .
\]
Moreover, the asymptotic formulae
\begin{equation}
\mu_n^{(s)} = \mu_n + o (s) \, , \quad v_n^{(s)} (t) = s \cos n t + o (s)
\label{s_a}
\end{equation}
hold for these solutions as $|s| \to 0$. Finally, each curve $C_n$ is of class
$C^1$.
\end{theorem}

\begin{figure}
\vspace{-4mm} 
\begin{center}
\SetLabels
 \L (0.5*0.01) $\mu$\\
 \L (0.56*0.41) $C_1$\\
 \L (-0.06*0.5) $\| v \|_\infty$\\
\endSetLabels 
\leavevmode \strut\AffixLabels{\includegraphics[width=68mm]{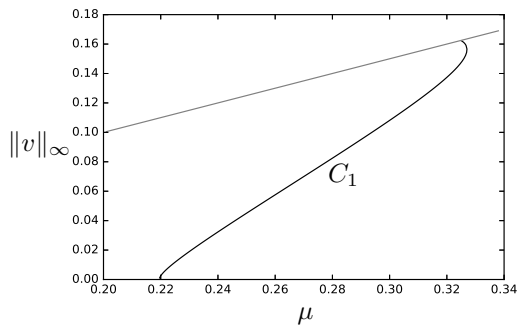}}
\end{center}
\vspace{-3mm} \caption{The branch of solutions of equation \eqref{37} with $r=4/5$,
bifurcating from the zero solution at $\mu_1 (4/5) = 0.219512195122$. The upper bound
mentioned prior to Definition~1 is also included.} \vspace{-4mm}
\label{fig:1}
\end{figure}

The last assertion is a consequence of the fact that ${\cal F}_{v v}$ is continuous
which is obvious. This theorem is illustrated in Fig.~1, where we have a plot of the
bifurcation branch $C_1$ in terms of $\mu$ and the norm of solution $\|v\|_\infty$
in the space $L^\infty (-\pi, \pi)$. The plotted branch bifurcating from $\mu_1 (4/5)$
has no secondary bifurcation points as the analogous branch for equation
\eqref{bid}; see \cite{BDT1,BDT2} for the rigorous proof and detailed discussion.
Moreover, it exhibits the phenomenon of a turning point at the largest value of
$\mu$ attained on $C_1$, occurring high on the branch; see further details in
Section~4.3. (The fastest traveling wave of given period corresponds to this point.)
By means of a different method this property was demonstrated by \cite{CN}, whereas
our method shows that it also takes place for equation \eqref{bid} on the branch
bifurcating from $\mu_1 (0)$. This phenomenon is related to the `Tanaka instability'
found numerically in \cite{Tan}, and later investigated analytically in \cite{Saf}.

\subsection{Solutions of Babenko's equation define periodic waves}

Let us outline a procedure demonstrating how to obtain a solution of problem
\eqref{lapp}--\eqref{bep} from that of Babenko's equation; that is, if equation
\eqref{37} with $r \in (0, 1)$ is satisfied by some $\mu > 0$ and an even function
$v$ (the existence of such pairs\,---\,at least in the form
\eqref{s_a}\,---\,follows from Theorem~2), then one can find the following:

(1) a $2 \, \pi$-periodic, symmetric curve with zero mean and a negative number
$-h$, which define the wave profile and the level of horizontal bottom,
respectively, thus giving a one-period water domain, say $\Omega$, on the $(x,
y)$-plane;

(2) a function $\psi$ harmonic in $\Omega$ and vanishing on its top side and two
positive constants serving as the right-hand side terms in the boundary conditions
\eqref{bcp} and \eqref{bep}.

Let we have an even, $2 \pi$-periodic solution $v$ of equation \eqref{37}, whose
Fourier coefficients we denote $b_0, b_1, b_2, \dots$ to distinguish these
coefficients from those in \eqref{z_u}, and let the periodic extension of $v$ to
$\RR$ be real-analytic. Using these coefficients, we define the following
holomorphic function on $A_r$:\\[-3mm]
\begin{equation}
z (u) = \ii \Big[ \log u - B + \sum_{k=1}^\infty b_k \left( u^k - r^{2 k} u^{-k}
\right) \Big] \, . \label{z_u'}
\end{equation}
Here $B$ is a real number that will be determined below in terms of the Fourier
coefficients of $v$. Let us consider the images that correspond under this mapping
to the curves and segments of $\partial A_r$. First, we see that $z (\E^{\ii t}) = x
(t) + \ii y (t)$ for $t \in (-\pi, \pi)$, where\\[-3mm]
\begin{equation}
x (t) = - t - \sum_{k=1}^\infty b_k \left( 1 + r^{2 k} \right) \sin kt \, , \ \ y
(t) = - B + \sum_{k=1}^\infty b_k \left( 1 - r^{2 k} \right) \cos kt \, .
\label{x_y}
\end{equation}
Since this curve given parametrically serves as the upper part of $\partial \Omega$,
we require its mean value to vanish. This gives that\\[-3mm]
\begin{equation}
B = \frac{1}{2} \sum_{k=1}^\infty k b_k^2 \left( 1 - r^{4 k} \right) \, ,
\label{d}
\end{equation}
where the series converges because the Fourier coefficients of the real-analytic $v$
decay faster than any power of $k$. 

Now we are in a position to determine the mean depth of flow $h$. In view of
symmetry we have that $z (r)$ is the mid-point of the bottom; that is, $z (r) = -
\ii h$. Then putting $u=r$ into \eqref{z_u'}, we find that\\[-3mm]
\begin{equation}
h = B - \log r = \frac{1}{2} \sum_{k=1}^\infty k b_k^2 \left( 1 - r^{4 k} \right) -
\log r \, , \label{h}
\end{equation}
and so is positive; here the last equality is a consequence of \eqref{d}. Thus, the
second expression \eqref{x_y} takes the following form:\\[-3mm]
\begin{equation}
y (t) = - (h + \log r) + \sum_{k=1}^\infty b_k \left( 1 - r^{2 k} \right) \cos kt \,
. \label{_y}
\end{equation}
Hence the curve $z_s = \{ x = x (t) = - t - (\mathcal{B}_r \, y) (t) , \ \ y = y
(t) ; \ \ t \in [-\pi, \pi] \}$ has the zero mean value. Here, the first formula
\eqref{HTB_r} is applied to express $x (t)$ in terms of $y (t)$.

Furthermore, we have\\[-3mm]
\begin{equation}
z (|u| \E^{\pm \ii \pi}) = \mp \pi + \ii \Big[ \log |u| - h + \sum_{k=1}^\infty
(-1)^k b_k \left( u^{2 k} - r^{2 k} \right) \! / |u|^{k} \Big] \quad \mbox{for} \ u
\in [-1, -r] \, , \label{right}
\end{equation}
thus obtaining two vertical segments $z_-$ and $z_+$ on the lines $x = -\pi$ and $x
= \pi$ respectively. 

Taking into account \eqref{d} and \eqref{h}, we see that\\[-3mm]
\begin{equation}
z (r \E^{\ii t}) = - \ii h - t - 2 \sum_{k=1}^\infty b_k r^k \sin k t \quad
\mbox{for} \ t \in [-\pi, \pi]  \label{horiz}
\end{equation}
on the inner circumference. This defines a horizontal segment $z_b$ on the line $y =
-h$.

\begin{figure}
\begin{center}
\leavevmode \strut\AffixLabels{\includegraphics[width=48mm]{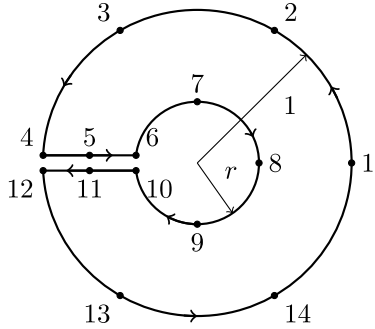}}
\end{center}
\vspace{-2mm} \caption{A sketch of the annular domain $A_r$ with several points on
its boundary marked in the counter-clockwise order.}
\label{fig:1}
\end{figure}

It is clear that the curve $\Gamma = z_+ \cup z_s \cup z_- \cup z_b$ constructed
above is closed and one can check (for example, numerically) that the set $\Omega$
enclosed within $\Gamma$ is a domain. The next step is to show that $z (u)$ defined
with the help of the Fourier coefficients of $v$ is a conformal mapping of $A_r$
onto $\Omega$. For this purpose one can use the boundary correspondence principle;
its form relevant for our case (see, for example, \cite{E}, Chapter~5, Theorem~1.3)
is formulated for the convenience of the reader.

\begin{theorem}
[The boundary correspondence principle] Let $D$ and $D^*$ be two bounded simply
connected domains with piecewise smooth boundaries and let $f$ be holomorphic in $D$
and continuous in $\bar D$. If $f (p)$ parametrises $\partial D^*$ counter-clockwise
provided $p$ is a counter-clockwise parametrisation of $\partial D$, then $f$ is a
conformal mapping of $D$ onto~$D^*$.
\end{theorem}

According to this theorem $z (u)$ maps $A_r$ onto $\Omega$ conformally provided one
can show (for example, numerically) that the map $\partial A_r \ni u \mapsto z \in
\Gamma$ is a homeomorphism. Moreover, condition \eqref{eta} is fulfilled for $\eta
(x) = y (t (x))$; here $t (x)$ is the inverse of $x = - t - (\mathcal{B}_r \, y)
(t)$, existing provided the curve $z_s$ is not self-intersecting. Thus, the curve $y
= \eta (x)$ defines the upper side of $\Omega$.

\begin{figure}
\begin{center}
\SetLabels
 \L (0.58*0.66) $\Gamma$\\
 \L (0.4*0.3) $\Omega$\\
\endSetLabels 
\leavevmode \strut\AffixLabels{\includegraphics[width=64mm]{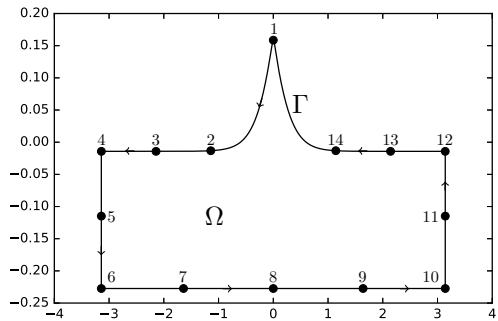}}
\end{center}
\vspace{-3mm} \caption{The curve $\Gamma$ corresponding to $\partial A_r$ through
the mapping $z (u)$ defined by \eqref{z_u'} and \eqref{d}, where the sequence
$\{b_k\}_{k=0}^\infty$ consists of the Fourier coefficients of $v$. The latter
solves \eqref{37} with $r=4/5$ and $\mu \approx 0.32671$, and $(\mu, v)$ belongs to
the branch bifurcating from $\mu_1$. The marked points on $\Gamma$ correspond to
those having the same numbers on $\partial A_r$ in Fig.~2. The mean depth of the
one-wave domain $\Omega$ is $h \approx 0.22739$, whereas the wave amplitude is
approximately equal to 0.17326.} \vspace{-6mm}
\label{fig:3}
\end{figure}

Figs~2--6 illustrate how the boundary correspondence principle works numerically in
recovering Stokes waves from solutions of \eqref{37}. We consider the equation with
$r=4/5$ and take the solution $(\mu, v)$ with $\mu \approx 0.32671$. This solution
belongs to the branch bifurcating from $\mu_1$ (equal to $0.219512195122$ for
$r=4/5$), and the value of $\mu$ in point is close to the critical one on this
branch (see Fig.~1 and Fig.~8). Substituting the Fourier coefficients of $v$ into
\eqref{z_u'} and \eqref{d}, we define $z (u)$ which is holomorphic in $A_r$ and maps
$\overline {A_r}$ onto $\overline \Omega$ continuously; the latter set is the
closure of the prospective one-wave domain. To demonstrate that $z (u)$ is a
conformal mapping we choose several points on $\partial A_r$, numbering them
counter-clockwise (see Fig.~2), and calculate their images on $\partial \Omega$,
assigning to each the same number as the object point has on $\partial A_r$. It
occurs that the images are also numbered counter-clockwise in agreement with the
boundary correspondence principle (see Fig.~3).

To be sure that the counter-clockwise boundary correspondence is not violated
between the chosen points we provide three figures more. In Fig.~4, the graph
of\\[-3mm]
\begin{equation}
x_h (t) = - t - 2 \sum_{k=1}^\infty b_k r^k \sin k t \label{horiz'}
\end{equation}
is plotted for $r=4/5$ and $t$ varying from 0 to $\pi$ (this parametrises the upper
half of the inner circumference clockwise provided it is considered as a part of
$\partial A_r$; see Fig.~2). According to \eqref{horiz}, this gives the left-hand
half of the bottom shown in Fig.~3 also parametrised clockwise. Since
\eqref{horiz'} is a monotonic function, there is no violation of the boundary
correspondence on the bottom because, by symmetry, it is sufficient to check this on
its right-hand half only.

\begin{figure}
\vspace{-4mm}
\begin{center}
\SetLabels
 \L (0.5*0.02) $t$\\
 \L (0.0*0.5) $x_h$\\
\endSetLabels 
\leavevmode \strut\AffixLabels{\includegraphics[width=68mm]{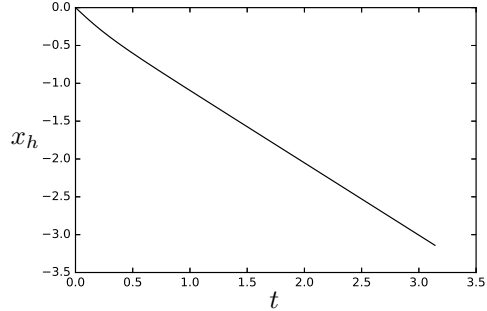}}
\end{center}
\vspace{-4mm} \caption{The graph of \eqref{horiz'} with $r=4/5$; its monotonicity
confirms that the boundary correspondence is not violated on the bottom part of
$\Gamma$.} \vspace{-4mm}
\label{fig:4}
\end{figure}

\begin{figure}
\begin{center}
\SetLabels
 \L (0.486*0.0) $|u|$\\
 \L (0.0*0.5) $y_+$\\
\endSetLabels 
\leavevmode \strut\AffixLabels{\includegraphics[width=68mm]{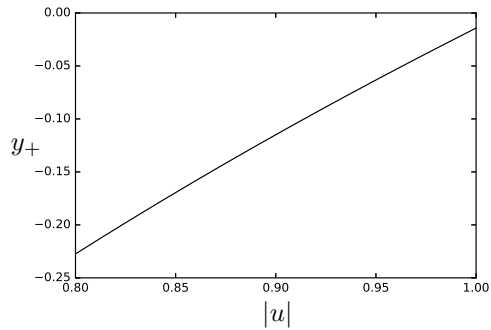}}
\end{center}
\vspace{-2mm} \caption{The graph of \eqref{right'} with $r=4/5$; its monotonicity
confirms that the boundary correspondence is not violated on the right-hand side of
$\Gamma$.} \vspace{-4mm}
\label{fig:5}
\end{figure}

In Fig.~5, the graph of
\begin{equation}
y_+ (u) = \log |u| - h + \sum_{k=1}^\infty (-1)^k b_k \left( u^{2 k} - r^{2 k}
\right) \! / |u|^{k} \label{right'}
\end{equation}
is plotted for $r=4/5$ and $u$ varying from $-r$ to $-1$ (this parametrises the lower
side of the cut counter-clockwise provided it is considered as a part of $\partial
A_r$; see Fig.~2). According to \eqref{right}, this gives the right-hand side of
$\Gamma$ shown in Fig.~3. Since \eqref{right'} is a monotonic function, there is no
violation of the boundary correspondence on the right-hand side of $\Gamma$.

\begin{figure}
\vspace{-4mm}
\begin{center}
\SetLabels
 \L (0.5*0.0) $t$\\
 \L (0.02*0.5) $x$\\
\endSetLabels 
\leavevmode \strut\AffixLabels{\includegraphics[width=68mm]{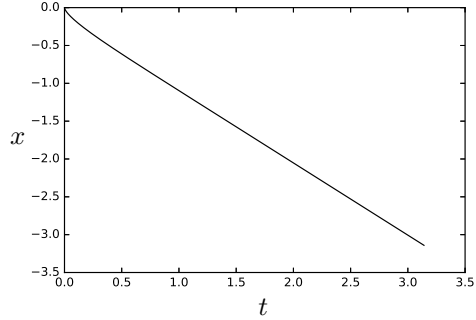}}
\end{center}
\vspace{-4mm}
\caption{The graph of the first function \eqref{x_y} with $r=4/5$; its monotonicity
confirms that the boundary correspondence is not violated on the left-hand half of
the upper part of $\Gamma$.} \vspace{-4mm}
\label{fig:6}
\end{figure}

Finally, the graph of the first function \eqref{x_y} is plotted in Fig.~6 for
$r=4/5$ and $t$ varying from $0$ to $\pi$ (this parametrises the upper half of the
exterior circumference of $\partial A_r$ counter-clockwise; see Fig.~2). According
to the first equation \eqref{x_y}, this parametrises the left-hand part of the upper
side of $\Gamma$ shown in Fig.~3. Since \eqref{x_y} is a monotonic function, there
is no violation of the boundary correspondence on this part of $\Gamma$.

It remains to check that $\Omega$ is a one-wave domain for some Stokes wave; that
is, there exists a stream function $\psi$ defined on $\overline \Omega$ so that it
satisfies conditions \eqref{bcp}--\eqref{bep} with some constant serving as the
right-hand side term in \eqref{bcp}, whereas $\mu$ stands in \eqref{bep}. For this
purpose we map $\Omega$ conformally on an auxiliary rectangle
\[ R^* = \{ (\varphi^*, \psi^*) : -\pi < \varphi^* < \pi , -\psi_0 < \psi^* < 0 \} 
\]
so that the images of $z_s$ and $z_b$ are the top and bottom parts of $\partial R^*$
respectively, whereas the value $\psi_0 > 0$ will be be chosen later. Thus, there
are harmonic functions $\varphi^*$ and $\psi^*$ defined on $\Omega$, and for every
$\psi_0$ the image of $R^*$ under the mapping $\E^{-\ii (\varphi^* + \ii \psi^*)}$
is the annular domain $A_\rho$ with some $\rho$. It is clear that the value of
$\rho$ decreases from unity to zero as $\psi_0$ characterising $R^*$ increases from
zero to infinity. Requiring $\rho$ to be equal to $r$, we fix the value of $\psi_0$,
thus determining $\psi^*$ which, in its turn, gives the constant value $-Q_*$ that
stands on the right-hand side of condition \eqref{bcp}; here the sign is chosen so
that $Q_*$ is positive. It should be noted that this procedure guarantees that
condition \eqref{kcp} is also fulfilled. It remains to use $\varphi^*$ and $\psi^*$
for determining $\psi$ so that it satisfies condition \eqref{bep} along with
\eqref{bcp} and \eqref{kcp}.

Using the Fourier coefficients $b_1, b_2, \dots$ of $v$ in formula \eqref{33}, we
obtain the function $z_{\varphi^*}$ holomorphic in $A_r$ and non-vanishing there.
According to equation \eqref{37}, we have that
\[ \left[ \{1 - 2 \mu^{-1} y (u)\} \overline{z_{\varphi^*} (u)} - 1 \right]_{|u|=1}
\]
is the limit as $|u| \to 1$ of some holomorphic function given in $A_r$ and having
its imaginary part equal to zero on $\partial A_r \cap \{ |u| = r \}$. Besides, the
same property holds for $z_{\varphi^*}$, and so it is also true for the function
whose limit as $|u| \to 1$ is equal to
\[ \left[ \{1 - 2 \mu^{-1} y (u)\} |z_{\varphi^*} (u)|^2 \right]_{|u|=1} \, .
\]
Therefore, we have that\\[-3mm]
\[ 1 - 2 \mu^{-1} \eta (x) = q^2 |\nabla \psi^* (x, \eta (x))|^2 , \quad x \in (-\pi,
\pi) ,
\]
with some $q > 0$ and $\eta$ defined above. For $\psi = q \sqrt \mu \, \psi^*$ the
last relation coincides with \eqref{bep}.

Thus, the triple $(\mu, \eta, \psi)$ satisfies problem P$(Q_0,h)$ with $h$ defined
by \eqref{h}, whereas $Q_0 = q \sqrt \mu \, Q_*$ and $Q_*$ depends on $r$ implicitly
as described above. This completes the description of a procedure how to obtain a
solution of problem \eqref{lapp}--\eqref{bep} from the given solution of Babenko's
equation.

\section{Numerical solution of Babenko's equation}

In this section, we describe a numerical method for solving equation \eqref{37} in
the class of even, periodic functions on $(-\pi, \pi)$. The existence of small
solutions of this kind is proved in Section~3.2, whereas general solutions are
discussed in Section~5. The essence of our method is to calculate the solution's
Fourier coefficients $b_0, b_1, \ldots$, which allows us to restore the conformal
mapping $z(u)$ (see Section~3.3), thus demonstrating numerically the equivalence of
Babenko's equation and problem P$(Q_0,h)$.

\subsection{Transformation of \eqref{37} to a form convenient for discretisation}
 
Let $r \in [0, 1)$ be fixed, then $J_r = \mathcal B_r \D / \D t$ is a self-adjoint
operator on $L^2_{per} (-\pi, \pi)$ of $2 \pi$-periodic square integrable functions.
Its domain is $H_0$ (see Section~3.2 for the definition), and it can also be defined
by linearity from $J_r \cos n t = \lambda_n \cos nt$ for $n=0,1,\dots$ and $J_r \sin
n t = \lambda_n \sin nt$ for $n=1,2,\dots$; here the eigenvalues are $\lambda_n =
\mu_n^{-1}$ for $n \geq 1$ and $\lambda_0 = 0$; cf. \eqref{lambda_n}. Since the
corresponding eigenfunctions form a basis in $L^2 (-\pi, \pi)$, the following
spectral decomposition holds:\\[-3mm]
\begin{equation}
J_r = \sum_{n = 1}^{\infty} \lambda_n ( \hat P_n + \tilde P_n ) .
\label{spect_decomp}
\end{equation}
Here $\hat P_n$ $(\tilde P_n)$ is the projector onto the subspace spanned by $\cos
nt$ ($\sin nt$, respectively).

Seeking solutions of \eqref{37} in $\hat H_0$, it is convenient to write the
equation in an equivalent form to accelerate numerical calculations. This form
arises after replacing $J_r$ in \eqref{37} by the right-hand side of
\eqref{spect_decomp} with omitted $\tilde P_n$, which is possible in view of the
bijection between $\hat H_0$ and the Sobolev space $W^{1,2} (0, \pi)$; indeed, for
every $w \in W^{1,2} (0, \pi)$ its even extension $v$ belongs to $\hat H_0$ and vice
versa. Therefore, it is convenient to put $\mathcal J_r = \sum_{n = 1}^{\infty}
\lambda_n P_n$, where $P_n$ is the projector onto the subspace of $L^2 (0, \pi)$
spanned by $\cos nt$. Then $\mathcal J_r$ is defined on $W^{1,2} (0, \pi)$ and
$\mathcal J_r w = J_r v (= \mathcal B_r v')$ almost everywhere on $(0, \pi)$, and so
\eqref{37} takes the following equivalent form\\[-3mm]
\begin{equation}
\mu \mathcal J_r w = w + w \mathcal J_r w + \frac 12 \mathcal J_r (w^2) \, , \quad t
\in (0, \pi) \, , \label{spectralBabenko}
\end{equation}
where $w (t)$ is sought in $W^{1,2} (0, \pi)$. To solve this equation numerically, a
modified version of the software SpecTraVVave is applicable; the latter is available
freely at the site indicated in \cite{MVK}, whereas its detailed description can be
found in \cite{KMV}.

For the reason made clear below, we amend \eqref{spectralBabenko} further; namely,
we set $\mu_0 = 1$ and put $\mathcal L_r = \sum_{n = 0}^{\infty} \mu_n P_n$. Hence
$\mathcal L_r$ is invertible and ${\mathcal L_r}^{-1} = P_0 + \mathcal J_r$; that
is, $\mathcal L_r \mathcal J_r = I - P_0$, where $I$ is the identity operator.
Applying $\mathcal L_r$ to both sides of \eqref{spectralBabenko}, we obtain the
following equation:
\begin{equation}
\mu (I - P_0) w = \mathcal L_r w + \mathcal L_r ( w \mathcal J_r w ) + \frac 12 (I -
P_0) w^2 \, , \quad t \in (0, \pi) \, . \label{inverse_spectralBabenko}
\end{equation}
It should be noted that the unbounded operator $\mathcal J_r$ is present in the
nonlinear part of the last equation only, and so one can expect that
\eqref{inverse_spectralBabenko} would demonstrate better numerical stability.
Finally, equations \eqref{inverse_spectralBabenko} and \eqref{37} are equivalent in
the following sense. The sets $\{ b_n (w) \}_{n=0}^\infty$ and $\{ b_n (v)
\}_{n=0}^\infty$ of the Fourier coefficients coincide for solutions
of\eqref{inverse_spectralBabenko} and \eqref{37}, respectively, provided the value
of $\mu$ is the same for both solutions.

For equation \eqref{inverse_spectralBabenko} the existence of small solutions
follows from its equivalence to \eqref{37}. It can also be established directly with
the help of the Crandall--Rabinowitz theorem; see Section~3.2, which yields the
asymptotic formulae \eqref{s_a} for the branch of solutions of
\eqref{inverse_spectralBabenko} bifurcating from $\mu_n$ and trivial $w$. This can
serve as a solution guess to start with in the numerical procedure.

\vspace{-1mm}

\subsection{Discretisation of equation \eqref{inverse_spectralBabenko}}

We use the standard cosine collocation method, according to which solutions of
\eqref{inverse_spectralBabenko} are are sought in the form of linear combinations of
$\cos mx$, $m = 0, 1, \dots$\,---\,a basis in $L^2(0, \pi)$. For the discretisation
the subspace $\mathcal S_N$ spanned by the first $N$ cosines is used, which is
defined by their values at the collocation points $x_n = \pi \frac{2n - 1}{2N}$ for
$n = 1, \ldots, N$. Thus, for any $f \in W^{1,2} (0, \pi)$ the vector $f^N$ given by
its coordinates\\[-2mm]
\[ f^N_n = \sum_{k=0}^{N-1} (P_k f) (x_n) \, , \quad n = 1, \dots, N ,
\]
is considered. The operator $\mathcal L_r^N$, discretising $\mathcal L_r$, is
defined as follows:\\[-2mm]
\[ ( \mathcal L_r^N f^N )_n = \sum_{k=0}^{N-1} (P_k \mathcal L_r f) (x_n) \, , \quad 
n = 1, \dots, N .
\]
Furthermore, $\mathcal J_r^N$ and $P_0^N$ are introduced as the discretisations of
$\mathcal J_r$ and $P_0$ respectively.

These definitions are correct because $f^N$ defines the function $f$ with values
$f(x_n) = f^N_n$ uniquely up to a projection on the subspace orthogonal to $\mathcal
S_N$. It is clear that each of these discrete operators is a composition of the
discrete cosine transform, some diagonal matrix and the inverse discrete cosine
transform. The diagonal matrix for $\mathcal L_r^N$ is $\{ 1, \ldots, \mu_{N-1} \}$,
whereas the diagonal for $\mathcal J_r^N$ is $\{ 0, \lambda_1, \ldots, \lambda_{N-1}
\}$, and $\{ 1, 0, \ldots, 0 \}$ is the diagonal for $P_0^N$. The discrete analogue
of \eqref{inverse_spectralBabenko} is as follows:\\[-3mm]
\begin{equation}
\mathcal L_r^N w^N - \mu \left( I - P_0^N \right) w^N + \mathcal L_r^N \left( w^N
\mathcal J_r^N w^N \right) + \frac 12 \left( I - P_0^N \right) \left( w^N \right)^2
= 0 . \label{discrete_Babenko}
\end{equation}
Since solutions $(\mu, w^N)$ of this equation form curves in the $(\mu, a)$-plane,
where\\[-3mm]
\[ a = \| w^N \| = \max_n |w^N_n| ,
\]
it is convenient to parametrise these curves for making calculations more effective.
Thus, due to a new parameter, say $\theta$, we have $\mu = \mu(\theta)$ and $a =
a(\theta)$ on each curve of solutions. Therefore, $\mu(\theta)$ must be substituted
into \eqref{discrete_Babenko} instead of $\mu$, and this algebraic system must be
complemented by the equation:\\[-4mm]
\begin{equation}
\max _{n = 1, \ldots, N} |w^N_n| = a(\theta) . \label{waveheight}
\end{equation}
The resulting system \eqref{discrete_Babenko}--\eqref{waveheight} has $N + 1$
equations with the following unknowns $\theta, w^N_1, \ldots, w^N_N$. Hence the
standard Newton's iteration method is applicable for finding bifurcations from a
trivial solution, and the Crandall--Rabinowitz asymptotic formula \eqref{s_a} yields
an initial guess. Further details concerning the proposed parametrisation and the
particular realisation of algorithm can be found in \cite{KMV}.

\begin{figure}
\vspace{-4mm}
\begin{center}
\SetLabels
 \L (0.5*0.01) $\mu$\\
 \L (0.764*0.41) $C_2$\\
 \L (0.81*0.8) \tiny $C_{21}$\\
 \L (0.428*0.31) $C_3$\\
 \L (0.43*0.56) \tiny $C_{31}$\\
 \L (0.256*0.26) $C_4$\\
 \L (0.232*0.43) \tiny $C_{41}$\\
 \L (-0.08*0.44) $\| v \|_\infty$\\
\endSetLabels 
\leavevmode
\strut\AffixLabels{\includegraphics[width=68mm]{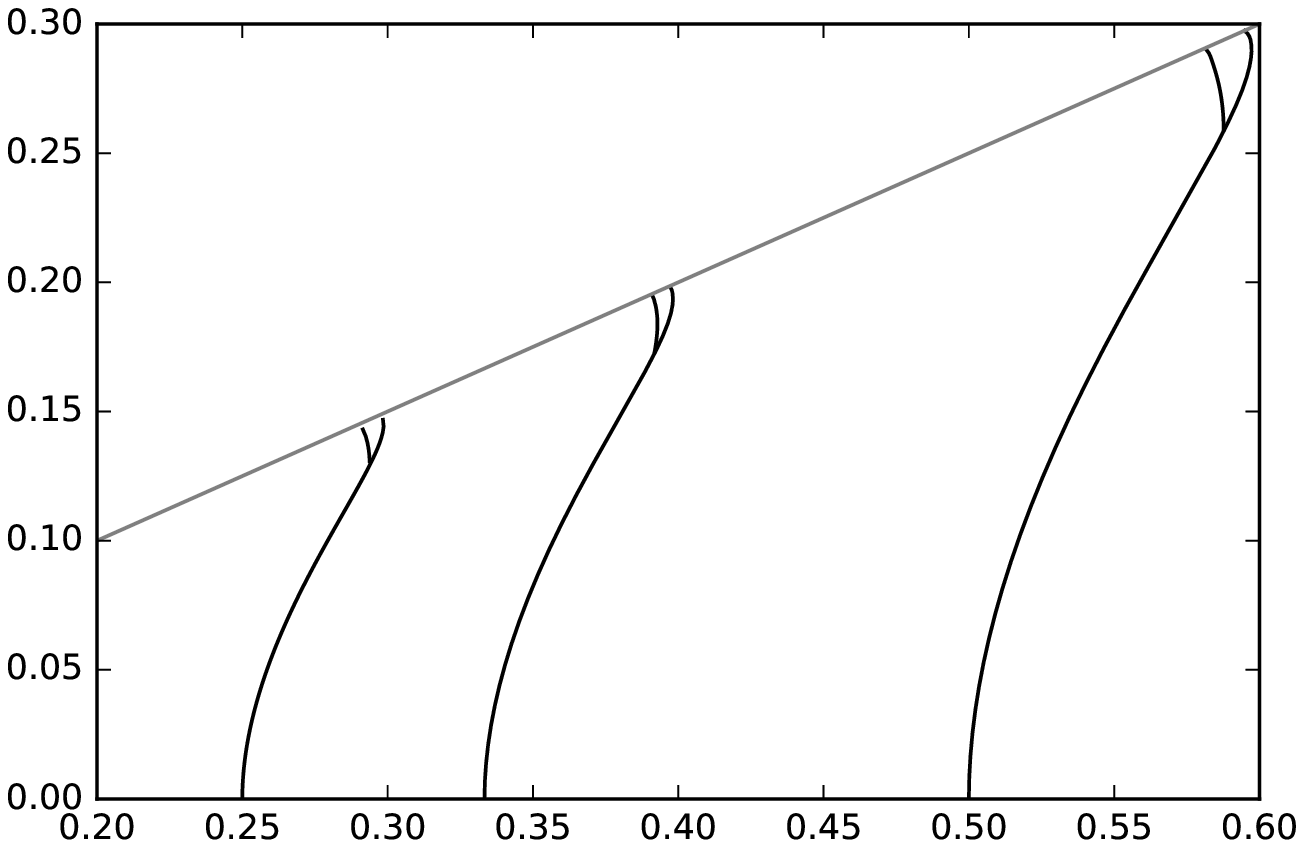}}
\end{center}
\vspace{-4mm}
\caption{The solution branches $C_2$, $C_3$ and $C_4$ for equation \eqref{37} with
$r=0$, bifurcating from the zero solution at $\mu_2 (0) = 1/2$, $\mu_3 (0) = 1/3$
and $\mu_4 (0) = 1/4$ respectively. The secondary solution branches are denoted
$C_{21}$, $C_{31}$ and $C_{41}$ respectively. The upper bound mentioned prior to
Definition~1 is also included.} \vspace{-3mm}
\label{fig:7}
\end{figure}

\begin{figure}
\begin{center}
\SetLabels
 \L (0.5*0.01) $\mu$\\
 \L (0.428*0.31) $C_1$\\
 \L (-0.08*0.44) $\| v \|_\infty$\\
\endSetLabels 
\leavevmode
\strut\AffixLabels{\includegraphics[width=68mm]{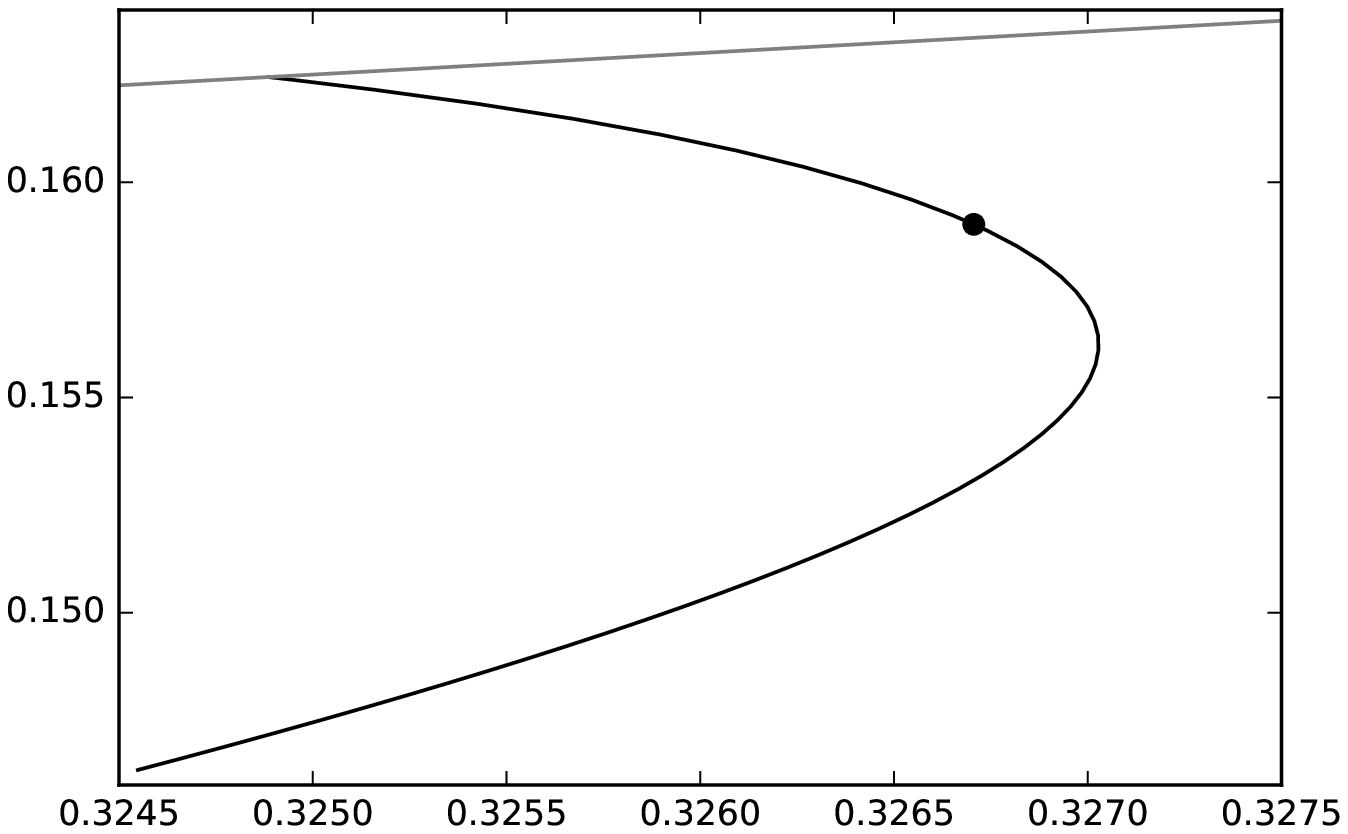}}
\end{center}
\vspace{-4mm} \caption{The solution branch $C_1$ for equation \eqref{37} with
$r=4/5$ in a vicinity of the turning point, whose characteristics are as follows:
$\mu \approx 0.32671$ the solution's $L^\infty$-norm is approximately equal to
$0.15862$. The bold dot marks the solution plotted in Fig.~3. The upper bound
mentioned prior to Definition~1 is also included.} \vspace{-3mm}
\label{fig:8}
\end{figure}

\subsection{Bifurcation curves for equation \eqref{37}}

We begin with the results of a test calculation in which the algorithm described in
Section~4.2 is applied to equation \eqref{inverse_spectralBabenko} with $r = 0$,
thus giving bifurcation curves for equation \eqref{bid}. The curves plotted in
Fig.~7 show the bifurcations from a trivial solution and the first three secondary
bifurcations for this case; the curve $C_1$ is omitted because its behaviour is
similar to that presented in Fig.~1, including the presence of a turning point. The
secondary bifurcation branches $C_{21}$, $C_{31}$ and $C_{41}$ bifurcate from $C_2$,
$C_3$ and $C_4$, respectively, at the points, where $\mu$ is approximately equal to
$0.58768$, $0.39172$ and $0.29389$ respectively. These values are in good agreement
with those presented by Aston \cite{A}; see Table~1 in his paper.

\begin{figure}
\vspace{-2mm}
\begin{center}
\SetLabels
 \L (0.5*0.01) $\mu$\\
 \L (0.5*0.41) $C_3$\\
 \L (0.74*0.764) \tiny $C_{31}$\\
 \L (-0.08*0.5) $\| v \|_\infty$\\
\endSetLabels 
\leavevmode \strut\AffixLabels{\includegraphics[width=68mm]{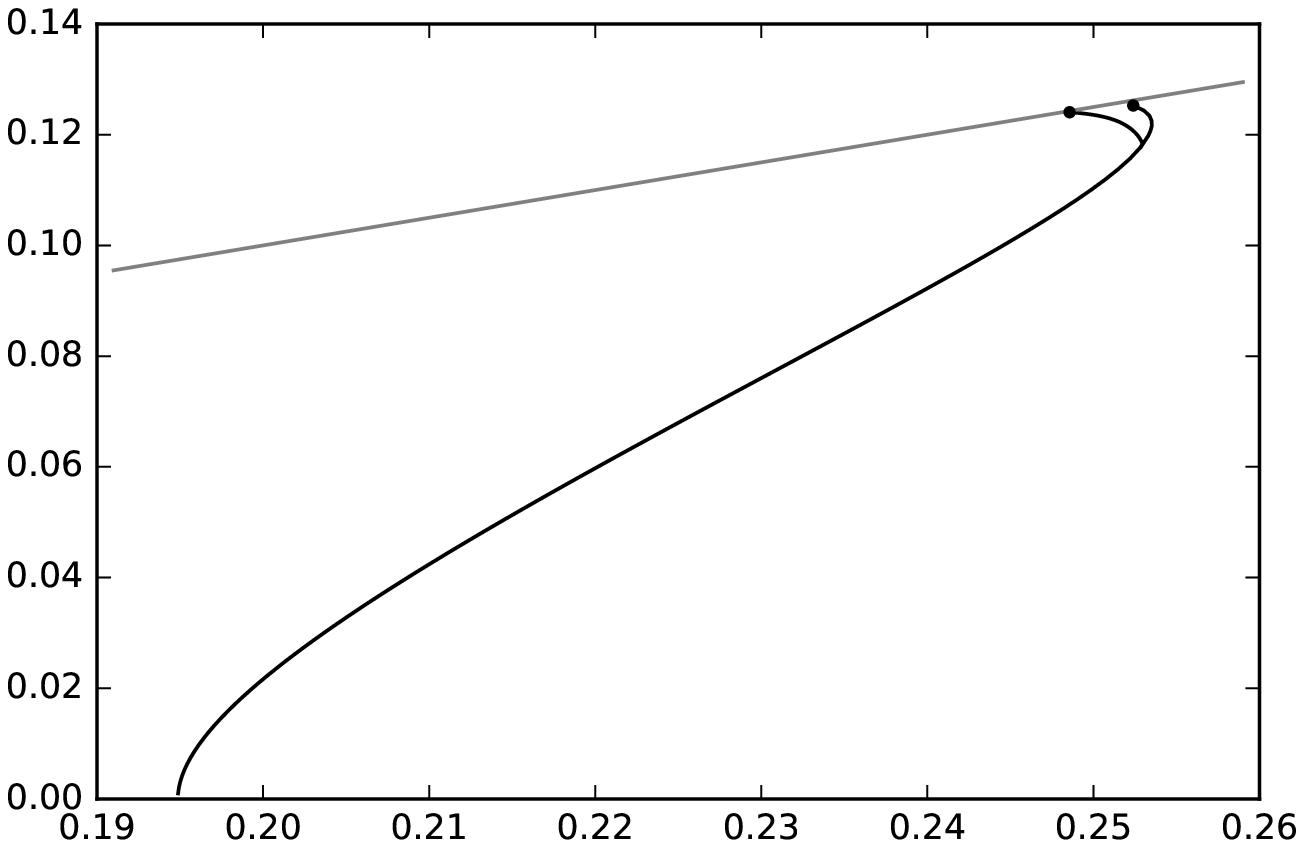}}
\end{center}
\vspace{-4mm}
\caption{The branch $C_3$ of solutions of equation \eqref{37} with $r=4/5$,
bifurcating from the zero solution at $\mu_3 (4/5) = 0.194868414381$. The secondary
solution branch $C_{31}$ bifurcates from $C_3$ at $\mu \approx 0.25298$. The dots
mark those solutions on $C_3$ and $C_{31}$, whose wave profiles are plotted in
Fig.~10 and Fig.~11 respectively. The upper bound mentioned prior to Definition~1 is
also included.} \vspace{-4mm}
\label{fig:9}
\end{figure}

Now we turn to numerical results obtained for equation \eqref{37} with $r=4/5$. The
so\-lution branch $C_1$ is presented in Fig.~1, and some of its characteristics are
described after Theorem~2. In particular, it is pointed out that it has a turning
point, and so we give a zoomed plot of the curve $C_1$ in a vicinity of this point;
see Fig.~8, where bold dot marks one of two solutions corresponding to $\mu \approx
0.32671$. The wave profile corresponding to this solution is plotted in Fig.~3,
where some of its characteristics are given; moreover, its $L^\infty$-norm is
approximately equal to $0.15862$.

The last example concerns the solution branch $C_3$ for equation \eqref{37} with
$r=4/5$. It is presented in Fig.~9, where one observes the presence of a turning
point as well as the secondary bifurcation. Indeed, the branch $C_{31}$ bifurcates
from $C_3$ at the point, where $\mu$ is approximately equal to $0.25298$, and
shortly after that $C_3$ has its turning point. The algorithm proposed in
Section~4.2 allows us to solve \eqref{inverse_spectralBabenko} up to both critical
values on $C_3$ and $C_{31}$; see Fig.~10 and Fig.~11, respectively, for the plots
of wave profiles corresponding to these solutions.

In Fig.~10, the wave profile corresponds to the end-point solution on the branch
$C_3$; $\mu \approx 0.25175$ for this solution of equation \eqref{37} with $r=4/5$.
Like a small-amplitude wave characterised by the second formula \eqref{s_a}, this
profile has the wavelength $2 \pi / 3$, and so three wave periods are plotted.
Moreover, this Stokes wave has the extreme form; that is, the tangents to two smooth
arcs form the angle $2 \pi / 3$ at every crest. The tangency is demonstrated with
sufficient accuracy in the figure, where the angle inscribed into the wave profile
has the sides $y = y_c \pm x / \sqrt 3$ with $y_c = y (0)$; see \eqref{_y} for $y
(t)$ and the first formula \eqref{x_y} for $x (t)$ that describe the profile
parametrically. Of course, the same tangency with similar angles takes place at
every crest.

\begin{figure}
\vspace{-4mm}
\begin{center}
\SetLabels
 \L (0.5*0.01) $x$\\
 \L (0.04*0.42) $y$\\
\endSetLabels 
\leavevmode
\strut\AffixLabels{\includegraphics[width=68mm]{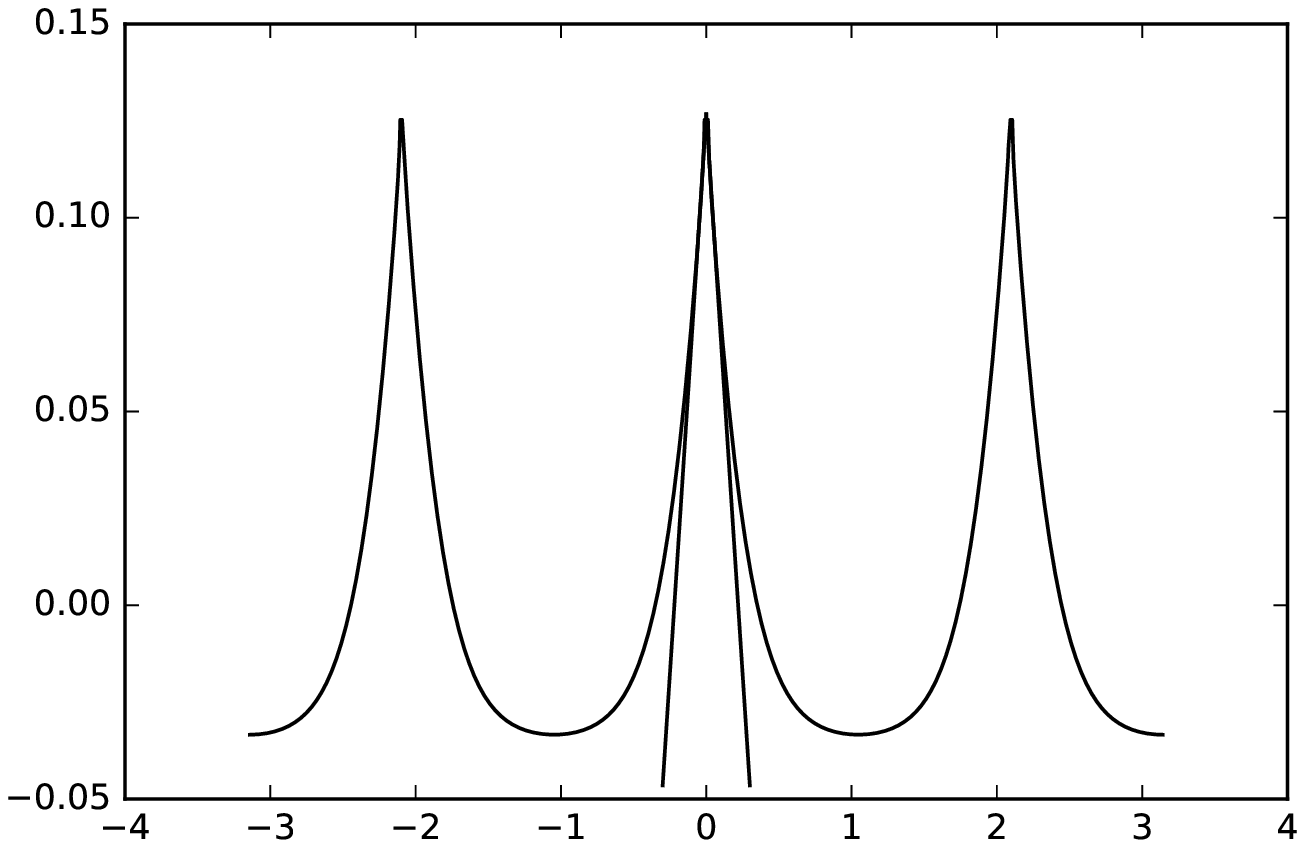}}
\end{center}
\vspace{-4mm}
\caption{The wave profile of the extreme form corresponding to the
end-point solution on the branch $C_3$ for equation \eqref{37} with $r=4/5$. The
characteristics of this wave are as follows: $\mu \approx 0.25175$; the profile's
crests (troughs) are at $y = y_c \approx 0.12777$ ($y = y_t \approx -0.03312$
respectively).} \vspace{-6mm}
\label{fig:10}
\end{figure}

\begin{figure}
\begin{center}
\SetLabels
 \L (0.5*0.01) $x$\\
 \L (0.0*0.46) $y$\\
\endSetLabels 
\leavevmode
\strut\AffixLabels{\includegraphics[width=68mm]{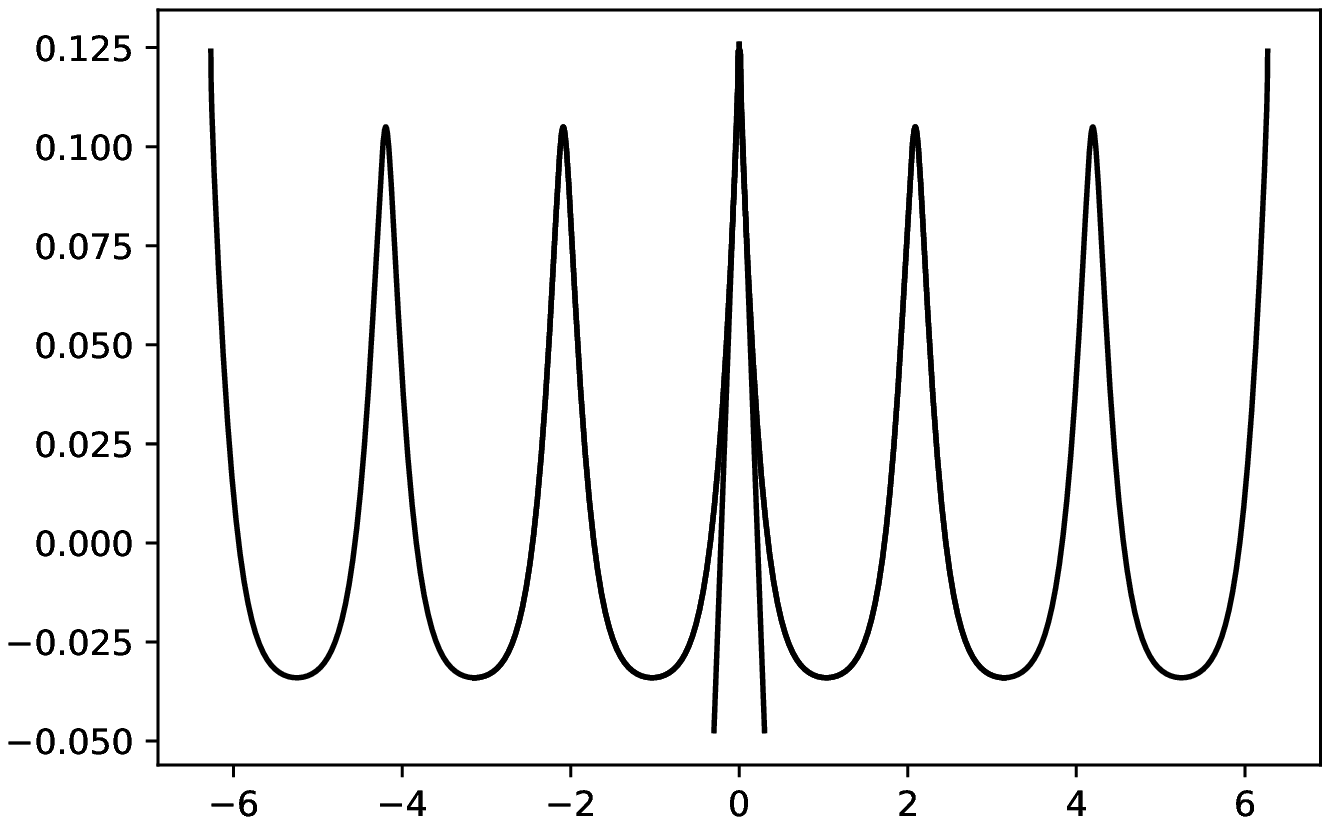}}
\end{center}
\vspace{-4mm}
\caption{The wave profile of the extreme form corresponds to the end-point solution
on the branch $C_{31}$ for equation \eqref{37} with $r=4/5$. Its characteristics are
as follows: $\mu \approx 0.24827$ the profile's smooth crests (troughs) are at $y =
\tilde y_c \approx 0.10406$ ($y = y_t \approx -0.03310$ respectively), whereas the
peaks are at $y = \hat y_c \approx 0.12608$.} \vspace{-6mm}
\label{fig:11}
\end{figure}

In Fig.~11, the wave profile corresponds to the end-point solution on the branch
$C_{31}$; $\mu \approx 0.24827$ for this solution of equation \eqref{37} with
$r=4/5$. The profile has the wavelength $2 \pi$, and so two wave periods are
plotted. Thus, the period-tripling occurs as $C_{31}$ bifurcates from the branch
$C_3$; an analogous effect is described in \cite{Zuf} for waves on infinitely deep
water (see, in particular, Fig.~3 on p.~25 of his paper). Moreover, the wave is
symmetric with respect to the vertical through the highest, mid-period crest. The
latter has the extreme form like every crest in Fig.~10, whereas the wave profile is
smooth at two other crests on the period.

\section{Concluding remarks}

We have considered the nonlinear problem describing steady, gravity waves on water
of finite depth. This problem is reduced to a single pseudo-differential operator
equation~\eqref{37} (Babenko's equation), which generalises the well-known equation
\eqref{bid} describing waves on infinitely deep water. The local bifurcation for
\eqref{37} is investigated analytically with the help of the Crandall--Rabinowitz
theorem, whereas a combination of analytical and numerical methods is applied for
demonstrating that the initial, free-boundary problem and Babenko's equation are
equivalent in the following sense. For every solution of the initial problem one of
its components, namely, the free-surface elevation, is a solution of Babenko's
equation for some value of the parameter on which the equation's operator depends;
this value is determined by the solution of the free-boundary problem. On the
contrary, every solution of Babenko's equation defines a solution of some
free-boundary problem through a certain procedure.

Besides, we outline an algorithm which allows us to solve Babenko's equation
numerically using a modification of the free software SpecTraVVave; see \cite{MVK}.
It should be emphasised that the developed numerical procedure is not only very
fast, but remarkable for its high accuracy. The latter is essential when computing
solutions to which wave profiles of the extreme form correspond, thus allowing us to
plot global bifurcation branches presented in Section~4.3.

This paper is just an initial step in studies of Babenko's equation both
analytically and numerically. First, it is desirable to prove rigorously that every
solution of Babenko's equation defines a solution of the free-boundary problem that
describes steady waves on a flow of finite depth with certain characteristics.
Second, it is natural to show that the profiles of waves below the highest, that has
the extreme form being non-smooth at its highest point, are real analytic curves.
Third, one has to demonstrate the absence of sub-harmonic bifurcations in a
neighbourhood of every point, where the bifurcation from the zero solution occurs.
Finally, a global Stokes-wave theory should be developed and used for proving that
there exist sub-harmonic bifurcations on branches of smooth waves close to the
highest wave. All these results had been established for waves on infinitely deep
water on the basis of equation \eqref{bid}; see \cite{BDT1,BDT2}.

An interesting direction for further numerical investigations is to find higher
bifurcations that might exist for waves on water of finite depth as it happens in
the case of deep water as had been shown in \cite{A}, where just several isolated
points of higher bifurcations are listed in Table~1. Since the algorithm based on
equation \eqref{37} and realised by using the software SpecTraVVave is a rather
robust tool, one could apply it for calculating branching bifurcation curves that
have more than one point of bifurcation.

In conclusion, we outline what equation \eqref{37} has in common with Babenko's
equation for finite depth obtained by Constantin, Strauss and V\u{a}rv\u{a}ruc\u{a}
\cite{CSV}; see Remark 4 in their paper. Namely,\\[-3mm]
\begin{equation}
\tilde \mu {\mathcal C}_d (\tilde v') = \tilde v + \tilde v {\mathcal C}_d (\tilde
v') + {\mathcal C}_d (\tilde v' \tilde v) \label{CSV}
\end{equation}
literally coincides with (2.50) in \cite{CSV} with one exception. We use $d$ as the
operator's parameter instead of $h$. There are two reasons for this: (1) $d$ and $h$
are equal to each other in Remark 4, since $k$ is taken equal to unity there; (2) a
different quantity is denoted by $h$ in Section~2.2 of our paper and $h$ will be
used in that meaning below.

It is obvious that the form of the last equation is exactly the same as that of
\eqref{37}, but what about the meaning of symbols involved? First, the parameter $d
> 0$ is equal to the so-called conformal mean depth and the latter is defined
uniquely by the water domain; see \cite{CV}, Appendix~A. However, this depth,
generally speaking, is not equal to the non-dimensional mean depth of the water
domain $D$ introduced in Section~2.2; see, in particular, formulae \eqref{dlv}. By
analogy with the conformal mean depth, it would be natural to call the parameter $r
\in (0, 1)$, on which the operator ${\mathcal B}_r$ depends in \eqref{37}, the
conformal mean radius of the water domain $D$. Furthermore, the conjugation operator
${\mathcal C}_d$ is defined for $2 \pi$-periodic functions on $\RR$ as follows. If
$f$ has zero mean value over a $2 \pi$ interval, that is, its Fourier series has the
form\\[-3mm]
\[ f (x) = \sum_{n=1}^\infty (a_n \cos n x + b_n \sin n x) , \quad x \in \RR ,
\]
then
\[ ({\mathcal C}_d f) (x) = \sum_{n=1}^\infty \coth n d \, (a_n \sin n x - b_n 
\cos n x) , \quad x \in \RR .
\]
This definition is similar to that of ${\mathcal B}_r$ in \eqref{HTB_r}, but with
the multiplier $\coth n d$ instead of $(1 + r^{2 n}) / (1 - r^{2 n})$. Moreover,
${\mathcal C}_d$ has the representation analogous to $\mathcal{B}_r = \mathcal{C} +
\mathcal{K}_r$ with $\mathcal{K}_r$ given by \eqref{K_r}; see formulae (A.9) and
(A.12) in \cite{CSV}. Thus, there is a significant similarity between ${\mathcal
C}_d$ and ${\mathcal B}_r$. The essential point that distinguishes ${\mathcal C}_d$
and $\mathcal{B}_r$ is that the latter operator is defined for all $2 \pi$-periodic
functions, whereas the domain of ${\mathcal C}_d$ is orthogonal to constants.

Finally, let us demonstrate that if $k = 1$ (this is the case in \cite{CSV}, Remark
4), the free surface profile satisfies the assumptions made in this paper (see
Section~2.1) and $\tilde v$ is an even and $2 \pi$-periodic solution of \eqref{CSV},
then $\tilde \mu = \mu$. Thus, the bifurcation parameter is the same in both
\eqref{CSV} and \eqref{37}.

According to Remark 4 in \cite{CSV}, we have\\[-3mm]
\begin{equation}
\tilde{\mu} = \frac{2R}{g} - 2 d  - 2 \beta , \label{14j}
\end{equation}
where $R$ is the Bernoulli constant in \eqref{4}, whereas the exact value of $\beta$
is unimportant for what follows. Since $k = 1$, formulae used in Section~2.2,
dealing with the derivation of the non-dimensional problem, imply that $l = \pi$ and
$H = h$, and so\\[-3mm]
\[ \mu = \frac{2R}{g} - 2 h .
\]
Combining this formula and \eqref{14j}, we see that $\tilde \mu = \mu$ holds, if we
show that $h = d + \beta$.

In order to prove the last equality, we first notice that the definition of
${\mathcal C}_d$ implies that ${\mathcal C}_d (\tilde v')$ and ${\mathcal C}_d
(\tilde v' \tilde v)$ are orthogonal to constants. Then equation \eqref{CSV} yields
that
\[ \int_{-\pi}^{\pi} \tilde v (x) [ 1 + {\mathcal C}_d (\tilde v') (x) ] \, \D x
= 0 \quad \Longleftrightarrow \quad \int_0^{\pi} \tilde v (x) [ 1 + {\mathcal C}_d
(\tilde v') (x) ] \, \D x = 0 \, ,
\]
where the second relation is a consequence of the assumption that $\tilde v$ is
even. To transform this relation we consider the parametric representation of the
free surface profile used in~\cite{CSV}:
\begin{equation}
X (x) = U (x, 0) = x + \mathcal{C}_d (\tilde{v} + \beta) \, , \quad Y (x) = V (x, 0)
= \tilde{v} + d + \beta \, , \quad x \in \RR , \label{XY}
\end{equation}
see (2.7), (2.8), (2.10) and Remark~4. Hence we have
\[ \int_0^{\pi} \tilde v (x) \frac{\D X}{\D x} (x) \, \D x = \int_0^{\pi} \tilde v 
(x (X)) \, \D X = 0 \, ,
\]
where it is taken into account that $X (x)$ is invertible on $(0, \pi)$ since $k =
1$. Averaging the second formula \eqref{XY} over $(0, \pi)$, we obtain
\[ h = \frac{1}{\pi} \int_0^\pi Y(x (X)) \, \D X = \frac{1}{\pi} \int_0^\pi \tilde{v}
(x (X)) \, \D X + d + \beta = d + \beta \, ,
\]
which yields the required equality $\tilde{\mu} = \mu$.

\vspace{1mm}

\noindent {\bf Acknowledgements.}

\noindent The authors are grateful to Henrik Kalisch without whose support the paper
would not appear. E.\,D. acknowledges the support from the Norwegian Research
Council.


\vspace{-2mm}

\begin{thebibliography}{99}

\bibitem{O} A. L. Afendikov, L. R. Volevich, G. P. Voskresenskii, I. M. Gelfand, A.
V. Zabrodin, O. V. Lokutsievskii, O. A. Oleinik, V. M. Tikhomirov, N. N. Chentsov,
Konstantin Ivanovich Babenko (obituary). {\sl Russian Math. Surveys} {\bf 43}
(1988), 139--151.

\bibitem{B} K. I. Babenko, Some remarks on the theory of surface waves of finite
amplitude. {\sl Soviet Math. Doklady} {\bf 35} (1987), 599--603.

\bibitem{BB} K. I. Babenko, A local existence theorem in the theory of surface waves
of finite amplitude. {\sl Soviet Math. Doklady} {\bf 35} (1987), 647--650.

\bibitem{BP} K. I. Babenko, V. Yu. Petrovich, A. I. Rakhmanov, A computational
experiment in the theory of surface waves of finite amplitude. {\sl Soviet Math.
Doklady} {\bf 38} (1989), 327--331.

\bibitem{BPR} K. I. Babenko, V. Yu. Petrovich, A. I. Rakhmanov, On a demonstrative
experiment in the theory of surface waves of finite amplitude. {\sl Soviet Math.
Doklady} {\bf 38} (1989), 626--630.

\bibitem{S} G. G. Stokes, On the theory of oscillatory waves. {\sl Camb. Phil. Soc.
Trans.} {\bf 8} (1847), 441--455.

\bibitem{PT} P. I. Plotnikov, J. F. Toland, Convexity of Stokes waves of extreme
form. {\sl Arch. Ration. Mech. Anal.} {\bf 171} (2004), 349--416.

\bibitem{N1} A. I. Nekrasov, On steady waves. {\sl Izvestia Ivanovo-Voznesensk.
Politekhn. Inst.} {\bf 3} (1921), 52--65; also {\sl Collected Papers, I}. Izdat.
Akad. Nauk SSSR, 1961, pp.~35--51 (both in Russian).

\bibitem{N3} A. I. Nekrasov, {\sl The Exact Theory of Steady Waves on the Surface of
a Heavy Fluid}. Izdat. Akad. Nauk SSSR, 1951; also {\sl Collected Papers, I}. Izdat.
Akad. Nauk SSSR, 1961, pp.~358--439 (both in Russian); translated as University of
Wisconsin MRC Report no.~813, 1967.

\bibitem{N2} A. I. Nekrasov, On steady waves on the surface of a heavy fluid. {\sl
Proc. All-Russian Congr. of Matematicians, Moscow}, (1928), 258--262 (in Russian).

\bibitem{AT2} C. J. Amick, J. F. Toland, On periodic water-waves and their
convergence to solitary waves in the long-wave limit. {\sl Phil. Trans. Roy. Soc.
Lond. A} {\bf 303} (1981), 633--669.

\bibitem{LC} T. Levi-Civita, D\'etermination rigoureuse des ondes permanentes
d'amplieur finie. {\sl Math. Ann.} {\bf 93} (1925), 264--314.

\bibitem{St} D. J. Struik, D\'etermination rigoureuse des ondes p\'eriodiques dans
un canal \`a profondeur finie. {\em Math. Ann.} {\bf 95} (1926), 595--634.

\bibitem{Z2} E. Zeidler, {\sl Nonlinear Functional Analysis and its Applications,
IV}. Springer-Verlag 1987.

\bibitem{Z1} E. Zeidler, {\sl Nonlinear Functional Analysis and its Applications,
I}. Springer-Verlag, 1985.

\bibitem{BT} B. Buffoni, J. F. Toland, {\sl Analytic Theory of Global Bifurcation:
an Introduction.} Princeton University Press, Princeton 2003.

\bibitem{T} J. F. Toland, Stokes waves. {\sl Topol. Methods Nonlinear Anal.} {\bf 7}
(1996), 1--48. Errata. Ibid {\bf 8} (1997), 412--413.

\bibitem{Kras} Yu. P. Krasovskii, On the theory of steady waves of finite amplitude.
{\sl USSR Comput. Math. Math. Phys.} {\bf 1} (1961), 996--1018.

\bibitem{BDT1} B. Buffoni, E. N. Dancer, J. F. Toland, The regularity and local
bifurcation of steady periodic waves. {\sl Arch. Ration. Mech. Anal.} {\bf 152}
(2000), 207--240.

\bibitem{BDT2} B. Buffoni, E. N. Dancer, J. F. Toland, The sub-harmonic bifurcation
of Stokes waves. {\em Arch. Ration. Mech. Anal.} {\bf 152} (2000), 241--271.

\bibitem{OS} H. Okamoto, M. Sh\={o}ji, {\sl The Mathematical Theory of Permanent
Progressive Water-Waves}. World Scientific, Singapore 2001.

\bibitem{Z} A Zygmund, {\sl Trigonometric Series, I \& II}. Cambridge University
Press, Cambridge 1959.

\bibitem{ST} E. Shargorodsky, J. F. Toland, Bernoulli free-boundary problems. {\sl
Memoirs AMS}\/ {\bf 96} (2008), no. 914.

\bibitem{LH1} M. S. Longuet-Higgins, Some new relations between Stokes’s
coefficients in the theory of gravity waves. {\sl J. Inst. Maths. Applics.} {\bf 22}
(1978), 261--273.

\bibitem{BS} J. G. B. Byatt-Smith, The equivalence of Bernoulli's equation and a set
of integral relations for periodic waves. {\sl IMA J. Appl. Math.} {\bf 23} (1979),
121--130.

\bibitem{LH2} M. S. Longuet-Higgins, Bifurcation in gravity
waves. {\sl J. Fluid Mech.} {\bf 151} (1985), 457--475.

\bibitem{Ba} A. M. Balk, A Lagrangian for water waves. {\sl Phys. Fluids} {\bf 8}
(1996), 416--420.

\bibitem{LH3} M. S. Longuet-Higgins, Lagrangian moments and mass transport in Stokes
waves. Part~2. Water of finite depth. {\sl J. Fluid. Mech.} {\bf 186} (1988),
321--336.

\bibitem{CS} B. Chen, P. G. Saffman, Numerical evidence for the existence of new
types of gravity waves of permanent form on deep water. {\sl Stud. Appl. Math.} {\bf
62} (1980), 1--21.

\bibitem{VB} J.-M. Vanden-Broeck, Some new gravity waves in water of finite depth.
{\sl Phys. Fluids} {\bf 26} (1983), 2385--2387.

\bibitem{VBS} J.-M. Vanden-Broeck, L. W. Schwartz, Numerical computation of steep
gravity waves in shallow water. {\sl Phys. Fluids} {\bf 22} (1979), 1868--1873.

\bibitem{CN} W. Craig, D. P. Nicholls, Travelling gravity water waves in two and
three dimensions. {\sl European J. Mech. B/Fluids} {\bf 21} (2002), 615--641.

\bibitem{Zuf1} J. A. Zufiria, Weakly nonlinear non-symmetric gravity waves on water
of finite depth. {\sl J. Fluid Mech.} {\bf 180} (1987), 371--385.

\bibitem{BM} C. Baesens, R. S. MacKay, Uniformly travelling water waves from a
dynamical systems viewpoint: some insights into bifurcations from Stokes' family.
{\sl J. Fluid Mech.} {\bf 241} (1992), 333--347.

\bibitem{CSV} A. Constantin, W. Strauss, E. V\u{a}rv\u{a}ruc\u{a}, Global
bifurcation of steady gravity water waves with critical layers. {\sl Acta Math.}
{\bf 217} (2016), 195--262.

\bibitem{Ben} T. B. Benjamin, Verification of the Benjamin--Lighthill conjecture
about steady water waves. {\sl J. Fluid Mech.} {\bf 295} (1995), 337--356.

\bibitem{KK} V. Kozlov, N. Kuznetsov, The Benjamin--Lighthill conjecture for steady
water waves (revisited). {\sl Arch. Ration. Mech. Anal.} {\bf 201} (2011), 631--645.

\bibitem{KK1} V. Kozlov, N. Kuznetsov, Bounds for arbitrary steady gravity waves on
water of finite depth. {\sl J. Math. Fluid Mech.} {\bf 11} (2009), 325--347.

\bibitem{KK2} V. Kozlov, N. Kuznetsov, Fundamental bounds for steady water waves.
{\sl Math. Ann.} {\bf 345} (2009), 643--655.

\bibitem{Bod} T. B. Bodnar', On steady periodic waves on the surface of a fluid of
finite depth. {\sl J.~Appl. Mech. Tech. Phys.} {\bf 52} (2011), 378--384.

\bibitem{CR} M. G. Crandall, P. H. Rabinowitz, Bifurcation from simple eigenvalues.
{\sl J. Func. Anal.} {\bf 8} (1971), 321--340.

\bibitem{Tan} M. Tanaka, The stability of steep gravity waves. {\sl J. Phys. Soc.
Japan} {\bf 52} (1983), 3047--3055.

\bibitem{Saf} P. G. Saffman, The superharmonic instability of finite amplitude water
waves. {\sl J. Fluid Mech.} {\bf 159} (1985), 169--174.

\bibitem{E} M. A. Evgrafov, {\sl Analytic functions.} Dover, New York 1978.

\bibitem{MVK} D. Moldabayev, O. Verdier, H. Kalisch, SpecTraVVave 2018. Free
software available at https://github.com/olivierverdier/SpecTraVVave

\bibitem{KMV} H. Kalisch, D. Moldabayev, O. Verdier, A numerical study of nonlinear
dispersive wave models with SpecTraVVave. {\sl Electronic J. Diff. Equations} {\bf
2017} (2017), 1--23.

\bibitem{A} P. J. Aston, Analysis and computation of symmetry-breaking bifurcation
and scaling laws using group theoretic methods. {\sl SIAM J. Math. Anal.} {\bf 22}
(1991), 181--212.

\bibitem{Zuf} J. A. Zufiria, Non-symmetric gravity waves on water of infinite depth.
{\sl J. Fluid Mech.} {\bf 181} (1987), 17--39.

\bibitem{CV} A. Constantin, E. V\u{a}rv\u{a}ruc\u{a}, Steady periodic water waves
with constant vorticity: regularity and local bifurcation. {\sl Arch. Ration. Mech.
Anal.} {\bf 199} (2011), 33--67.

\end{thebibliography}
\end{document}